\documentclass{article}
\usepackage{amssymb,amsthm,amsmath}
\begin{document}
\title{Strong Duality Theorem for Continuous-Time Linear Programming Problems}
\author{\vspace{2mm}Hsien-Chung Wu
\thanks{e-mail:\vspace{5mm}
hcwu@nknucc.nknu.edu.tw}\\
Department of Mathematics, \\
National Kaohsiung Normal University, Kaohsiung 802, Taiwan}
\date{}
\maketitle
\newtheorem{Thm}{Theorem}[section]
\newtheorem{Def}{Definition}[section]
\newtheorem{Lem}{Lemma}[section]
\newtheorem{Pro}{Proposition}[section]
\newtheorem{Rem}{Remark}[section]
\newtheorem{Cor}{Corollary}[section]
\newtheorem{Ex}{Example}[section]

\newenvironment{Proof}{\noindent {\bf Proof}.}
{\hspace{5mm}\rule{2.5mm}{2mm}\medskip\par}

\begin{abstract}

This paper is aimed to prove the strong duality theorem for continuous-time
linear programming problems in which the coefficients are assumed to be
piecewise continuous functions.
The previous paper proved the strong duality theorem for the case of
piecewise continuous functions in which the discontinuities are
the left-continuities. In this paper, we propose the completely different type of
discretized primal and dual problems that can be used to prove the strong duality
theorem for the general situation of discontinuities.

\vspace{3mm}

\vspace{3mm}

\end{abstract}

\section{Introduction}

In Wu \cite{wu,wu2}, the coefficients appeared in the discretized primal and
dual problems are the function values of the coefficient functions taken at
the right-end points of the subdivided intervals. This simple type of formulation
can just be used to prove the strong duality theorem for the case of
piecewise continuous functions in which the discontinuities are the left-continuities.
In this paper, we shall extend to prove the strong duality theorem for the general situation
of discontinuities. We shall propose the completely different type of formulation for the
discretized primal and dual problems. In this paper, the coefficients in the
discretized primal and dual problems will consider the infimum and supremum
of the coefficient functions on the subdivided intervals, which is more complicated
than that of considering the function values of the coefficient functions taken at
the right-end points of the subdivided intervals in Wu \cite{wu,wu2}.

The theory of continuous-time linear programming problem has received
considerable attention for a long time. Tyndall \cite{tyn65,tyn67}
treated rigorously a continuous-time linear programming problem with
constant matrices, which had originated from the ``bottleneck problem''
proposed by Bellman \cite{bel57}. Levinson \cite{lev} generalized the results
of Tyndall by considering time-dependent matrices in which the functions
appearing in the objective and constraints were assumed to be continuous on the time
interval $[0,T]$. Meidan and Perold \cite{mei}, Papageorgiou \cite{pap} and
Schechter \cite{sch} have also obtained
some interesting results for the continuous-time linear programming problem.
Anderson {\em et al.} \cite{and83,and94,and96}, Fleischer and Sethuraman \cite{fle}
and Pullan \cite{pul93,pul95,pul96,pul00,pul02} investigated a subclass of
continuous-time linear programming problems, which is called separated
continuous-time linear programming problem and can be used to model the
job-shop scheduling problems. Weiss \cite{wei} proposed a simplex-like algorithm to solve
the separated continuous-time linear programming problem.

This paper is organized as follows. In Section 2, the problem formulation is
presented, and the weak duality theorem is proved.
In Section 3, in order to study the strong duality theorem, we propose a perturbed
continuous-time linear programming problem. Many useful results that will be used
to prove the strong duality theorem are derived.
In Section 4, discretized problems are formulated in which the
partition of the time interval $[0,T]$ is not taken as equally dividing $[0,T]$.
We also derive many useful results that will be used
to prove the strong duality theorem.
In Section 5, the strong duality theorem is proved.

\section{Formulation}

Let $A$ be a matrix with entries denoted by $a_{ij}$. We define
$\parallel A\parallel =\sum_{i,j} |a_{ij}|$.
Let $L^{\infty}_{p}[0,T]$ be the space of all measurable and
essentially bounded functions from the compact interval $[0,T]$ into the
$p$-dimensional Euclidean space $\mathbb{R}^{p}$.
If $p=1$, Then, we simply write $L^{\infty}[0,T]$.
For $f\in L^{\infty}[0,T]$, we define
\[\parallel f\parallel_{\infty}=\mbox{ess}\sup_{t\in [0,T]}
|f(t)|=\inf\left\{k:|f(t)|\leq k\mbox{ a.e. in }[0,T]\right\},\]
where the Lebesgue measure is considered.
Therefore, we have $|f(t)|\leq\parallel f\parallel_{\infty}$
a.e. in $[0,T]$. For ${\bf f}=(f_{1},\cdots ,f_{p})\in L^{\infty}_{p}[0,T]$, we define
\[\parallel {\bf f}\parallel_{\infty}^{p}
=\max_{i=1,\cdots ,p}\parallel f_{i}\parallel_{\infty}.\]
We consider the following assumptions:
\begin{itemize}
\item ${\bf a}\in L^{\infty}_{q}[0,T]$ and ${\bf c}\in L^{\infty}_{p}[0,T]$;

\item $B$ and $K$ are time-dependent $p\times q$ matrices
defined on $[0,T]$ and $[0,T]\times [0,T]$, respectively, such that
each entry is in the spaces $L^{\infty}[0,T]$ and $L^{\infty}([0,T]\times [0,T])$,
respectively.
\end{itemize}
The continuous-time linear programming problem is formulated as follows:
\begin{eqnarray*}
(\mbox{CLP}^{*}) & \max & \int_{0}^{T} {\bf a}^{\top}(t){\bf z}(t)dt\\
& \mbox{subject to} & B(t){\bf z}(t)\leq {\bf c}(t)+\int_{0}^{t}
K(t,s){\bf z}(s)ds\mbox{ for all $t\in [0,T]$},\\
&& {\bf z}\in L^{\infty}_{q}[0,T]\mbox{ and }
{\bf z}(t)\geq {\bf 0}\mbox{ for all $t\in [0,T]$}.
\end{eqnarray*}
The dual problem of $(\mbox{CLP}^{*})$ is defined as follows:
\begin{eqnarray*}
(\mbox{DCLP}^{*}) & \min & \int_{0}^{T} {\bf c}^{\top}(t){\bf w}(t)dt\\
& \mbox{subject to} & B^{\top}(t){\bf w}(t)\geq {\bf a}(t)+\int_{t}^{T}
K^{\top}(s,t){\bf w}(s)ds\mbox{ for all $t\in [0,T]$},\\
&& {\bf w}\in L^{\infty}_{p}[0,T]\mbox{ and }
{\bf w}(t)\geq {\bf 0}\mbox{ for all $t\in [0,T]$}.
\end{eqnarray*}
In this paper, we shall consider the following problems:
\begin{eqnarray*}
\mbox{(CLP)} & \max & \int_{0}^{T} {\bf a}^{\top}(t){\bf z}(t)dt\\
& \mbox{subject to} & B(t){\bf z}(t)\leq {\bf c}(t)+\int_{0}^{t}
K(t,s){\bf z}(s)ds\mbox{ a.e. in $[0,T]$,}\\
&& {\bf z}\in L^{\infty}_{q}[0,T]\mbox{ and }
{\bf z}(t)\geq {\bf 0}\mbox{ a.e. in $[0,T]$}
\end{eqnarray*}
and
\begin{eqnarray*}
\mbox{(DCLP)} & \min & \int_{0}^{T} {\bf c}^{\top}(t){\bf w}(t)dt\\
& \mbox{subject to} & B^{\top}(t){\bf w}(t)\geq {\bf a}(t)+\int_{t}^{T}
K^{\top}(s,t){\bf w}(s)ds\mbox{ a.e. in $[0,T]$,}\\
&& {\bf w}\in L^{\infty}_{p}[0,T]\mbox{ and }
{\bf w}(t)\geq {\bf 0}\mbox{ a.e. in $[0,T]$},
\end{eqnarray*}
where the constraints are assumed to be satisfied in the sense of a.e. in $[0,T]$.
The weak duality theorem can be similarly established
although the primal and dual problems (CLP) and (DCLP)
are defined in the sense of a.e. in $[0,T]$.

\begin{Thm}{\label{optt198}}
{\em (Weak Duality Theorem)}
If ${\bf z}$ and ${\bf w}$ are any arbitrary feasible solutions of
the primal and dual problems {\em (CLP)} and {\em (DCLP)}, respectively, then
\[\int_{0}^{T} {\bf a}^{\top}(t){\bf z}(t)dt\leq\int_{0}^{T}
{\bf c}^{\top}(t){\bf w}(t)dt.\]
\end{Thm}
\begin{Proof}
According to the constrains of problems (CLP) and (DCLP), we have
\begin{equation}{\label{*clpeq56}}
\sum_{j=1}^{q}\int_{0}^{T}z_{j}(t)\left [a_{j}(t)-
\sum_{i=1}^{p}B_{ij}(t)w_{i}(t)+
\sum_{i=1}^{p}\int_{t}^{T} K_{ij}(s,t)w_{i}(s)ds\right ]dt\leq 0
\end{equation}
and
\begin{equation}{\label{*clpeq57}}
\sum_{i=1}^{p}\int_{0}^{T}w_{i}(t)\left [c_{i}(t)-\sum_{j=1}^{q}B_{ij}(t)z_{j}(t)
+\sum_{j=1}^{q}\int_{0}^{t} K_{ij}(t,s)z_{j}(s)ds\right ]dt\geq 0.
\end{equation}
By Fubini's theorem, we also have
\begin{equation}{\label{dclp112}}
\int_{0}^{T}\int_{t}^{T}K_{ij}(s,t)z_{j}(t)w_{i}(s)dsdt=
\int_{0}^{T}\int_{0}^{t}K_{ij}(t,s)z_{j}(s)w_{i}(t)dsdt.
\end{equation}
Therefore, in the vectorial form, we obtain
\begin{align*}
0 & \geq\int_{0}^{T}{\bf z}^{\top}(t)\left [{\bf a}(t)
-B^{\top}(t){\bf w}(t)+\int_{t}^{T}K^{\top}(s,t){\bf w}(s)ds\right ]dt\\
& \quad -\int_{0}^{T}{\bf w}^{\top}(t)\left [{\bf c}(t)-B(t){\bf z}(t)
+\int_{0}^{t}K(t,s){\bf z}(s)ds\right ]dt
\mbox{ (by (\ref{*clpeq56}) and (\ref{*clpeq57}))}\\
& =\int_{0}^{T}{\bf z}^{\top}(t){\bf a}(t)dt+\int_{0}^{T}
{\bf w}^{\top}(t)\left [-B(t){\bf z}(t)+\int_{0}^{t}K(t,s){\bf z}(s)ds
\right ]dt\mbox{ (by (\ref{dclp112}))}\\
& \quad -\int_{0}^{T}{\bf c}^{\top}(t){\bf w}(t)dt-\int_{0}^{T}
{\bf w}^{\top}(t)\left [-B(t){\bf z}(t)+\int_{0}^{t}K(t,s){\bf z}(s)ds\right ]dt\\
& =\int_{0}^{T}{\bf z}^{\top}(t){\bf a}(t)dt
-\int_{0}^{T}{\bf c}^{\top}(t){\bf w}(t)dt.
\end{align*}
This completes the proof.
\end{Proof}

In the sequel, we are going to prove the strong duality theorem between
(CLP) and (DCLP) although these problems are considered in the sense of a.e. in $[0,T]$.

\section{Perturbed Formulation}

Given any $\epsilon\geq 0$, we consider the following perturbed problems:
\begin{eqnarray*}
(\mbox{CLP}_{\epsilon}) & \max & \int_{0}^{T} {\bf a}^{\top}(t){\bf z}(t)dt\\
& \mbox{subject to} & B(t){\bf z}(t)\leq
\left [{\bf c}(t)+\mbox{\boldmath $\epsilon$}\right ]
+\int_{0}^{t} K(t,s){\bf z}(s)ds\mbox{ a.e. in $[0,T]$},\\
&& {\bf z}\in L^{\infty}_{q}[0,T]\mbox{ and }
{\bf z}(t)\geq {\bf 0}\mbox{ a.e. in $[0,T]$}
\end{eqnarray*}
and
\begin{eqnarray*}
(\mbox{DCLP}_{\epsilon}) & \min & \int_{0}^{T} {\bf c}^{\top}(t){\bf w}(t)dt\\
& \mbox{subject to} & B^{\top}(t){\bf w}(t)
\geq \left [{\bf a}(t)-\mbox{\boldmath $\epsilon$}\right ]
+\int_{t}^{T} K^{\top}(s,t){\bf w}(s)ds\mbox{ a.e. in $[0,T]$},\\
&& {\bf w}\in L^{\infty}_{p}[0,T]\mbox{ and }
{\bf w}(t)\geq {\bf 0}\mbox{ a.e. in $[0,T]$},
\end{eqnarray*}
where $\mbox{\boldmath $\epsilon$}$ is a vector with all entries $\epsilon$.
Although the vectors $\mbox{\boldmath $\epsilon$}$ in problems
$(\mbox{CLP}_{\epsilon})$ and $(\mbox{DCLP}_{\epsilon})$ have the different dimensions,
we use the same notation for convenience.
If the constraints of $(\mbox{CLP}_{\epsilon})$ and
$(\mbox{DCLP}_{\epsilon})$ are assumed to be satisfied for all $t\in [0,T]$,
then the corresponding problems are denoted by $(\mbox{CLP}_{\epsilon}^{*})$ and
$(\mbox{DCLP}_{\epsilon}^{*})$.

Since each entry of ${\bf a}$, ${\bf c}$, $B$ and $K$ is measurable and essentially bounded
in $[0,T]$ and $[0,T]\times [0,T]$, respectively, we define
\begin{eqnarray}
&& \tau =\max_{j=1,\cdots ,q}\parallel a_{j}\parallel_{\infty};
\mbox{ that is, }|a_{j}(t)|\leq\tau
\mbox{ a.e. in $[0,T]$},\label{clpeq8}\\
&& \zeta =\max_{i=1,\cdots ,p}\parallel c_{i}\parallel_{\infty};
\mbox{ that is, }|c_{i}(t)|\leq\zeta\mbox{ a.e. in $[0,T]$},\label{clpeq2}\\
&& \eta =\max_{i=1,\cdots ,p;j=1,\cdots ,q}\parallel K_{ij}\parallel_{\infty};
\mbox{ that is, }|K_{ij}(t,s)|\leq\eta\mbox{ a.e. in $[0,T]\times [0,T]$},\label{*clpeq2}\\
&& \nu =\max_{j=1,\cdots ,q}\sum_{i=1}^{p}\parallel K_{ij}\parallel_{\infty};
\mbox{ that is, }\sum_{i=1}^{p}|K_{ij}(t,s)|
\leq\nu\mbox{ a.e. in $[0,T]\times [0,T]$},\label{clpeq50}\\
&& \phi =\max_{i=1,\cdots ,p}\sum_{j=1}^{q}\parallel K_{ij}\parallel_{\infty};
\mbox{ that is, }\sum_{j=1}^{q}|K_{ij}(t,s)|
\leq\phi\mbox{ a.e. in $[0,T]\times [0,T]$}.\label{clpeq3}
\end{eqnarray}
Let $\{f_{k}\}_{k=1}^{\infty}$ be a sequence of functions in $L^{\infty}[0,T]$.
We say that the sequence $\{f_{k}\}_{k=1}^{\infty}$ is uniformly essentially bounded
in $[0,T]$ if and only if there exists a positive constant $C$ such that
$\parallel f_{k}\parallel_{\infty}\leq C$ for each $k$.
If $\{{\bf f}_{k}\}_{k=1}^{\infty}$ is a sequence of vector-valued functions,
Then, we say that the sequence $\{{\bf f}_{k}\}_{k=1}^{\infty}$ is uniformly
essentially bounded if and only if there exists a positive constant $C$ such that
$\parallel f_{ik}\parallel_{\infty}\leq C$ for each $i$ and $k$, where $f_{ik}$ is
the $i$th entry of ${\bf f}_{k}$. For $f\in L^{2}[0,T]$, we recall
\[\parallel f\parallel_{2}=\left (\int_{0}^{T}f^{2}(t)dt\right )^{1/2}.\]
Then, the sequence $\{f_{k}\}_{k=1}^{\infty}$ of real-valued functions is uniformly
bounded in $[0,T]$ with respect to $\parallel\cdot\parallel_{2}$ if and only if
there exists a positive constant $C$ such that
$\parallel f_{k}\parallel_{2}\leq C$ for each $k$.
The concept of uniform boundedness of the sequence of vector-valued
functions $\{{\bf f}_{k}\}_{k=1}^{\infty}$ can be similarly defined.
We also see that if the sequence $\{{\bf f}_{k}\}_{k=1}^{\infty}$ is uniformly
essentially bounded in $[0,T]$, then it is also uniformly
bounded in $[0,T]$ with respect to $\parallel\cdot\parallel_{2}$.

We denote by ${\cal Z}_{\epsilon}$ and ${\cal W}_{\epsilon}$
the feasible sets of problems $(\mbox{CLP}_{\epsilon})$ and
$(\mbox{DCLP}_{\epsilon})$, respectively.
We say that the feasible set ${\cal Z}_{\epsilon}$ of
$(\mbox{CLP}_{\epsilon})$ is uniformly
essentially bounded if and only if there exists a positive constant $C$
such that each feasible solution of $(\mbox{CLP}_{\epsilon})$ is essentially
bounded by $C$. We are going to provide the sufficient conditions to guarantee
that the feasible set ${\cal Z}_{\epsilon}$ of $(\mbox{CLP}_{\epsilon})$ is
uniformly essentially bounded.
Gronwall's lemma was provided by Levinson \cite{lev}. We can similarly prove it in the
sense of a.e. in $[0,T]$.

\begin{Lem}{\label{optl196}}
{\em (Gronwall's lemma)}
Suppose that the real-valued function $g$ is integrable in $[0,T]$ and
$g(t)\geq 0$ a.e. in $[0,T]$ $($resp. for all $t\in [0,T])$.
If there exist constants $\theta_{1}\geq 0$ and $\theta_{2}>0$ such that
\begin{equation}{\label{opteq134}}
g(t)\leq\theta_{1}+\theta_{2}\cdot\int_{0}^{t}g(s)ds\mbox{ a.e. in $[0,T]$
$($resp. for all $t\in [0,T])$},
\end{equation}
then $g(t)\leq\theta_{1}\cdot e^{\theta_{2}t}$ a.e in $[0,T]$
$($resp. for all $t\in [0,T])$.
\end{Lem}
\begin{Proof}
We are going to prove the case of a.e. on $[0,T]$.
For $t\in [0,T]$, we define
\[G(t)=\int_{0}^{t}g(s)ds.\]
Then, we see that $G$ is continuous on $[0,T]$ and
$G'(t)=g(t)$ a.e. on $[0,T]$ by Royden \cite{roy}.
From (\ref{opteq134}), we also have
\begin{equation}{\label{opteq135}}
g(t)\leq\theta_{1}+\theta_{2}G(t)\mbox{ a.e. on $[0,T]$}.
\end{equation}
Using (\ref{opteq135}), we also have
\begin{align*}
\frac{d}{dt}\left (e^{-\theta_{2}t}G(t)\right )
& =-\theta_{2}e^{-\theta_{2}t}G(t)+e^{-\theta_{2}t}g(t)\\
& \leq -\theta_{2}e^{-\theta_{2}t}G(t)+e^{-\theta_{2}t}\cdot\left (
\theta_{1}+\theta_{2}G(t)\right )
=\theta_{1}e^{-\theta_{2}t}\mbox{ a.e. on $[0,T]$}.
\end{align*}
By taking integration, for $t\in [0,T]$, we have
\begin{equation}{\label{opteq137}}
\int_{0}^{t}\frac{d}{dt}\left (e^{-\theta_{2}s}G(s)\right )ds\leq
\int_{0}^{t}\theta_{1}e^{-\theta_{2}s}ds.
\end{equation}
Since $e^{-\theta_{2}s}G(s)$ is continuous on $[0,T]$,
the Lebesgue integral and Riemann integral are identical
as given in (\ref{opteq137}). Therefore, we have
\[e^{-\theta_{2}t}G(t)-G(0)\leq-\frac{\theta_{1}}
{\theta_{2}}e^{-\theta_{2}t}+\frac{\theta_{1}}{\theta_{2}}\]
for each $t\in [0,T]$. Since $G(0)=0$, we have
\begin{equation}{\label{opteq136}}
G(t)\leq\frac{\theta_{1}}{\theta_{2}}\left (e^{\theta_{2}t}-1\right )
\end{equation}
for each $t\in [0,T]$. Using (\ref{opteq135}) and (\ref{opteq136}), we obtain
\[g(t)\leq\theta_{1}+\theta_{2}G(t)\leq\theta_{1}+\theta_{2}\cdot
\frac{\theta_{1}}{\theta_{2}}\left (e^{\theta_{2}t}-1\right )
=\theta_{1}e^{\theta_{2}t}\mbox{ a.e. on $[0,T]$}.\]
This completes the proof.
\end{Proof}

\begin{Pro}{\label{optt120*}}
Suppose that there exist real-valued functions $\lambda_{i}$ satisfying
$0\leq\lambda_{i}(t)\leq 1$ a.e. in $[0,T]$ $($resp. for all $t\in [0,T])$ for
$i=1,\cdots ,p$ and a constant $\sigma >0$ satisfying
\[\min_{j=1,\cdots ,q}\left\{\sum_{i=1}^{p}\lambda_{i}(t)B_{ij}(t)
\right\}\geq\sigma\mbox{ a.e. in $[0,T]$ $($resp. for all $t\in [0,T])$}.\]
If the problem $(\mbox{\em CLP}_{\epsilon})$ is feasible, then
each feasible solution ${\bf z}^{(\epsilon)}(t)$ is bounded satisfying
\begin{align}
\left |z_{j}^{(\epsilon)}(t)\right | & \leq\parallel {\bf z}^{(\epsilon)}(t)\parallel
\leq\frac{p\cdot (\zeta +\epsilon )}{\sigma}\cdot
\exp\left (\frac{p\cdot\phi\cdot t}{\sigma}\right )\nonumber\\
& \leq\frac{p\cdot (\zeta +\epsilon )}{\sigma}\cdot
\exp\left (\frac{p\cdot\phi\cdot T}{\sigma}\right )\mbox{ a.e. in $[0,T]$
$($resp. for all $t\in [0,T])$}\label{clpeq337}
\end{align}
for $j=1,\cdots ,q$, where the constants $\zeta$ and $\phi$ are given in
$(\ref{clpeq2})$ and $(\ref{*clpeq2})$, respectively.
In other words, the feasible set ${\cal Z}_{\epsilon}$ of
$(\mbox{\em CLP}_{\epsilon})$ is uniformly essentially bounded in $[0,T]$
when $\epsilon$ is fixed.
\end{Pro}
\begin{Proof}
We are going to prove the case of a.e. in $[0,T]$.
Let ${\bf z}^{(\epsilon)}$ be a feasible solution of
$(\mbox{CLP}_{\epsilon})$. According to the constraints, we have
\begin{align*}
& \sum_{j=1}^{q}\sum_{i=1}^{p}\lambda_{i}(t)\cdot B_{ij}(t)\cdot z_{j}^{(\epsilon)}(t)\\
& \quad\leq\sum_{i=1}^{p}\lambda_{i}(t)\cdot\left |c_{i}(t)+\epsilon\right |+
\int_{0}^{t}\sum_{j=1}^{q}\sum_{i=1}^{p}\lambda_{i}(t)\cdot
K_{ij}(t,s)\cdot z_{j}^{(\epsilon)}(s)ds\mbox{ a.e. in $[0,T]$}.
\end{align*}
Therefore, we obtain
\begin{align*}
\sigma\cdot\parallel {\bf z}^{(\epsilon)}(t)\parallel
& \leq\sum_{j=1}^{q}\left [\left |z_{j}^{(\epsilon)}(t)\right |\cdot
\sum_{i=1}^{p}\lambda_{i}(t)B_{ij}(t)\right ]
=\sum_{j=1}^{q}\sum_{i=1}^{p}\lambda_{i}(t)\cdot B_{ij}(t)\cdot
\left |z_{j}^{(\epsilon)}(t)\right |\\
& \leq\sum_{i=1}^{p}\lambda_{i}(t)\cdot\left |c_{i}(t)+\epsilon\right |+
\int_{0}^{t}\sum_{i=1}^{p}\left [\sum_{j=1}^{q}\lambda_{i}(t)\cdot
K_{ij}(t,s)\cdot\left |z_{j}^{(\epsilon)}(s)\right |ds\right ]\\
& \leq\sum_{i=1}^{p}\left |c_{i}(t)+\epsilon\right |+\int_{0}^{t}\sum_{i=1}^{p}
\left [\sum_{j=1}^{q}K_{ij}(t,s)\cdot\left |z_{j}^{(\epsilon)}(s)\right |ds\right ]\\
& \leq p\cdot (\zeta +\epsilon )+p\cdot\phi\cdot\int_{0}^{t}
\parallel {\bf z}^{(\epsilon)}(s)\parallel ds\mbox{ a.e. in $[0,T]$}.
\end{align*}
By Gronwall's Lemma~\ref{optl196}, we obtain
\[\parallel {\bf z}^{(\epsilon)}(t)\parallel\leq
\frac{p\cdot (\zeta +\epsilon )}{\sigma}\cdot
\exp\left (\frac{p\cdot\phi\cdot t}{\sigma}\right )
\leq\frac{p\cdot (\zeta +\epsilon )}{\sigma}\cdot
\exp\left (\frac{p\cdot\phi\cdot T}{\sigma}\right )
\mbox{ a.e. in $[0,T]$}.\]
This completes the proof.
\end{Proof}

The following lemmas are very useful.

\begin{Lem}{\label{optl138}}
{\em (Riesz and Sz.-Nagy \cite[p.64]{rie})}
Let $\{f_{k}\}_{k=1}^{\infty}$ be a sequence in $L^{2}[0,T]$.
If the sequence $\{f_{k}\}_{k=1}^{\infty}$ is uniformly bounded with respect to
$\parallel\cdot\parallel_{2}$, then there exists a subsequence
$\{f_{k_{r}}\}_{r=1}^{\infty}$ which weakly converges to some $f_{0}\in L^{2}[0,T]$.
In other words, for any $g\in L^{2}[0,T]$, we have
\[\lim_{r\rightarrow\infty}\int_{0}^{T}
f_{k_{r}}(t)g(t)dt=\int_{0}^{T}f_{0}(t)g(t)dt.\]
\end{Lem}

\begin{Lem}{\label{optl157}}
{\em (Levinson \cite{lev})}
If the sequence $\{f_{k}\}_{k=1}^{\infty}$ is uniformly bounded in $[0,T]$
with respect to $\parallel\cdot\parallel_{2}$ and weakly converges
to some $f_{0}\in L^{2}[0,T]$, then
\[f_{0}(t)\leq\limsup_{k\rightarrow\infty}f_{k}(t)\mbox{ a.e. in $[0,T]$}\]
and
\[f_{0}(t)\geq\liminf_{k\rightarrow\infty}f_{k}(t)\mbox{ a.e. in $[0,T]$}.\]
\end{Lem}

\begin{Pro}{\label{*clpp47}}
Consider the sequence $\{\epsilon_{k}\}_{k=1}^{\infty}$ with
$\epsilon_{k}\rightarrow 0+$ as $k\rightarrow\infty$. Assume that $B(t)\geq {\bf 0}$
a.e. in $[0,T]$. The following statements hold true.
\begin{enumerate}
\item [{\em (i)}] Suppose that each problem $(\mbox{\em CLP}_{\epsilon_{k}})$
is feasible, and that ${\bf z}^{(\epsilon_{k})}$ is a feasible solution of
problem $(\mbox{\em CLP}_{\epsilon_{k}})$ such that the sequence
$\{{\bf z}^{(\epsilon_{k})}\}_{k=1}^{\infty}$ is uniformly essentially bounded.
Then, there exists a subsequence $\{{\bf z}^{(\epsilon_{k_{r}})}\}_{r=1}^{\infty}$
which weakly converges to some feasible solution ${\bf z}^{(0)}\in L_{q}^{2}[0,T]$ of
$(\mbox{\em CLP}_{0})=\mbox{\em (CLP)}$. Moreover, there exists a feasible solution
$\bar{\bf z}$ of {\em (CLP)} such that $\bar{\bf z}(t)\geq {\bf 0}$ for all
$t\in [0,T]$ and $\bar{\bf z}(t)={\bf z}^{(0)}(t)$ a.e. in $[0,T]$.

\item [{\em (ii)}] Suppose that ${\bf c}(t)\geq {\bf 0}$ for all $t\in [0,T]$,
and that $K(t_{0},s)\geq {\bf 0}$ a.e. in $[0,T]$ for each fixed $t_{0}\in [0,T]$.
If the sequence $\{{\bf z}^{(\epsilon_{k})}\}_{k=1}^{\infty}$ of feasible solutions of
problem $(\mbox{\em CLP}_{\epsilon_{k}})$ is uniformly essentially bounded,
then there exists a subsequence $\{{\bf z}^{(\epsilon_{k_{r}})}\}_{r=1}^{\infty}$
which weakly converges to some feasible solution ${\bf z}^{(0)}\in L_{q}^{2}[0,T]$ of
$\mbox{\em (CLP)}$. Moreover, there exists a feasible solution $\bar{\bf z}$
of $(\mbox{\em CLP}^{*})$ such that $\bar{\bf z}(t)={\bf z}^{(0)}(t)$ a.e. in $[0,T]$.
\end{enumerate}
\end{Pro}
\begin{Proof}
To prove part (i), since the sequence $\{{\bf z}^{(\epsilon_{k})}\}_{k=1}^{\infty}$
is uniformly essentially bounded in $[0,T]$, it follows
that this sequence is also uniformly bounded in $[0,T]$ with respect to
$\parallel\cdot\parallel_{2}$. Let $z_{j}^{(\epsilon_{k})}$ be the $j$th entry of
${\bf z}^{(\epsilon_{k})}$. Using Lemma~\ref{optl138}, there exists a subsequence of
$\{z_{j}^{(\epsilon_{k})}\}_{k=1}^{\infty}$
which weakly converges to some $z_{j}^{(0)}(t)\in L^{2}[0,T]$.
Therefore, we can construct a vector-valued subsequence
$\{{\bf z}^{(\epsilon_{k_{r}})}\}_{r=1}^{\infty}$ of
$\{{\bf z}^{(\epsilon_{k})}\}_{k=1}^{\infty}$
such that $\{z_{j}^{(\epsilon_{k_{r}})}\}_{r=1}^{\infty}$ weakly converges to
$z_{j}^{(0)}$ for $j=1,\cdots ,q$. For each $i=1,\cdots ,p$, the constraints say that
\begin{equation}{\label{*clpeq59}}
\sum_{j=1}^{q}B_{ij}(t)\cdot z_{j}^{(\epsilon_{k_{r}})}(t)\leq c_{i}(t)
+\epsilon_{k_{r}}+\sum_{j=1}^{q}\int_{0}^{t} K_{ij}(t,s)\cdot
z_{j}^{(\epsilon_{k_{r}})}(s)ds\mbox{ a.e. in $[0,T]$}.
\end{equation}
Using Lemma~\ref{optl157}, we also have
\begin{equation}{\label{clpeq1}}
\limsup_{r\rightarrow\infty}z_{j}^{(\epsilon_{k_{r}})}(t)\geq
z_{j}^{(0)}(t)\geq\liminf_{r\rightarrow\infty}z_{j}^{(\epsilon_{k_{r}})}(t)\geq 0
\mbox{ a.e. in $[0,T]$.}
\end{equation}
Since $\epsilon_{k_{r}}\rightarrow 0$ as $r\rightarrow\infty$ and
$B(t)\geq {\bf 0}$ a.e. in $[0,T]$, from (\ref{*clpeq59}) and (\ref{clpeq1}),
by taking the limit superior and using the weak convergence, we obtain
\begin{equation}{\label{opteq148}}
B(t){\bf z}^{(0)}(t)\leq\limsup_{r\rightarrow\infty}
B(t){\bf z}^{(\epsilon_{k_{r}})}(t)
\leq {\bf c}(t)+\int_{0}^{t} K(t,s){\bf z}^{(0)}(s)ds
\mbox{ a.e. in $[0,T]$.}
\end{equation}
This shows that ${\bf z}^{(0)}$ is a feasible solution of (CLP).
Let $N_{0j}=\{t\in [0,T]:z_{j}^{(0)}(t)<0\}$ and $N_{0}=\bigcup_{j=1}^{q}N_{0j}$.
Let $N_{1}$ be the subset of $[0,T]$ such that the inequality
(\ref{opteq148}) is violated. We define $N=N_{0}\cup N_{1}$.
Then from (\ref{clpeq1}) and (\ref{opteq148}), we see that the set $N$ has measure zero.
Now, we define
\begin{equation}{\label{*clpeq51}}
\bar{\bf z}(t)=\left\{\begin{array}{ll}
{\bf z}^{(0)}(t) & \mbox{if $t\not\in N$}\\
{\bf 0} & \mbox{if $t\in N$}.
\end{array}\right .
\end{equation}
Then, we see that $\bar{\bf z}(t)\geq {\bf 0}$ for all $t\in [0,T]$
and $\bar{\bf z}(t)={\bf z}^{(0)}(t)$ a.e. in $[0,T]$.
For $t\not\in N$, from (\ref{opteq148}), we have
\begin{equation}{\label{*clpeq52}}
B(t)\bar{\bf z}(t)=B(t){\bf z}^{(0)}(t)
\leq {\bf c}(t)+\int_{0}^{t} K(t,s){\bf z}^{(0)}(s)ds
={\bf c}(t)+\int_{0}^{t} K(t,s)\bar{\bf z}(s)ds.
\end{equation}
This shows that $\bar{\bf z}$ is a feasible solution of (CLP).

To prove part (ii), under the assumptions of ${\bf c}(t)$ and $K(t,s)$,
it is obvious that the problem $(\mbox{CLP}_{\epsilon_{k}})$ is feasible
for each $\epsilon_{k}$ with the trivial feasible solution
${\bf z}(t)={\bf 0}$ for all $t\in [0,T]$.
We consider $\bar{\bf z}$ defined in (\ref{*clpeq51}).
For $t\in N$, we have $B(t)\bar{\bf z}(t)={\bf 0}$.
Since ${\bf z}^{(0)}(t)\geq {\bf 0}$ a.e. in $[0,T]$ and
$K(t_{0},s)\geq {\bf 0}$ a.e. in $[0,T]$ for each fixed $t_{0}\in [0,T]$,
we obtain
\[B(t)\bar{\bf z}(t)={\bf 0}\leq {\bf c}(t)+\int_{0}^{t} K(t,s){\bf z}^{(0)}(s)ds
={\bf c}(t)+\int_{0}^{t} K(t,s)\bar{\bf z}(s)ds.\]
By referring to (\ref{*clpeq52}) for $t\not\in N$,
we see that $\bar{\bf z}(t)$ satisfies all the constraints of
primal problem (CLP) for all $t\in [0,T]$. This completes the proof.
\end{Proof}

\begin{Pro}{\label{*clpr62}}
Assume that $B(t)\geq {\bf 0}$ a.e. in $[0,T]$. The following statements hold true.
\begin{enumerate}
\item [{\em (i)}] Suppose that the problem $(\mbox{\em CLP}_{\epsilon})$ is feasible.
For any uniformly essentially bounded sequence
$\{{\bf z}^{(k)}\}_{k=1}^{\infty}$ of feasible solutions of
$(\mbox{\em CLP}_{\epsilon})$, there exists a subsequence
$\{{\bf z}^{(k_{r})}\}_{r=1}^{\infty}$ which weakly converges to some feasible solution
${\bf z}^{(\epsilon)}\in L_{q}^{2}[0,T]$ of $(\mbox{\em CLP}_{\epsilon})$.
Moreover, there exists a feasible solution $\bar{\bf z}^{(\epsilon)}$ of
$(\mbox{\em CLP}_{\epsilon})$
such that $\bar{\bf z}^{(\epsilon)}(t)\geq {\bf 0}$ for all $t\in [0,T]$
and $\bar{\bf z}^{(\epsilon)}(t)={\bf z}^{(\epsilon)}(t)$ a.e. in $[0,T]$.

\item [{\em (ii)}] Suppose that ${\bf c}(t)\geq {\bf 0}$ for all $t\in [0,T]$,
and that $K(t_{0},s)\geq {\bf 0}$ a.e. in $[0,T]$ for each fixed $t_{0}\in [0,T]$.
Given any uniformly essentially bounded sequence
$\{{\bf z}^{(k)}\}_{k=1}^{\infty}$ of
feasible solutions of $(\mbox{\em CLP}_{\epsilon})$, there exists a subsequence
$\{{\bf z}^{(k_{r})}\}_{r=1}^{\infty}$ which weakly converges
to some feasible solution ${\bf z}^{(\epsilon)}\in L_{q}^{2}[0,T]$ of
$(\mbox{\em CLP}_{\epsilon})$. Moreover, there exists a feasible solution
$\bar{\bf z}^{(\epsilon)}$ of $(\mbox{\em CLP}_{\epsilon}^{*})$ such that
$\bar{\bf z}^{(\epsilon)}(t)={\bf z}^{(\epsilon)}(t)$ a.e. in $[0,T]$.
\end{enumerate}
\end{Pro}
\begin{Proof}
To prove part (i), since the sequence $\{{\bf z}^{(k)}\}_{k=1}^{\infty}$
is uniformly essentially bounded in $[0,T]$,
we see that this sequence is also uniformly bounded in $[0,T]$ with respect to
$\parallel\cdot\parallel_{2}$. Using Lemma~\ref{optl138}, there exists a subsequence of
$\{z_{j}^{(k)}\}_{k=1}^{\infty}$ which weakly converges to some $z_{j}^{(\epsilon)}\in L^{2}[0,T]$.
Therefore, we can construct a vector-valued subsequence
$\{{\bf z}^{(k_{r})}\}_{r=1}^{\infty}$ of $\{{\bf z}^{(k)}\}_{k=1}^{\infty}$
such that $\{z_{j}^{(k_{r})}\}_{r=1}^{\infty}$ weakly converges to
$z_{j}^{(\epsilon)}$ for $j=1,\cdots ,q$. For each $i=1,\cdots ,p$, the feasibility says that
\begin{equation}{\label{*clpeq59}}
\sum_{j=1}^{q}B_{ij}(t)\cdot z_{j}^{(k_{r})}(t)\leq c_{i}(t)+\epsilon
+\sum_{j=1}^{q}\int_{0}^{t} K_{ij}(t,s)\cdot z_{j}^{(k_{r})}(s)ds\mbox{ a.e. in $[0,T]$}.
\end{equation}
Using Lemma~\ref{optl157}, we also have
\begin{equation}{\label{clpeq1}}
\limsup_{r\rightarrow\infty}z_{j}^{(k_{r})}(t)\geq
z_{j}^{(\epsilon)}(t)\geq\liminf_{r\rightarrow\infty}z_{j}^{(k_{r})}(t)\geq 0
\mbox{ a.e. in $[0,T]$.}
\end{equation}
Since $B(t)\geq {\bf 0}$ a.e. in $[0,T]$, from (\ref{*clpeq59}) and (\ref{clpeq1}),
by taking the limit superior and using the weak convergence, we obtain
\begin{equation}{\label{opteq148}}
B(t){\bf z}^{(\epsilon)}(t)\leq\limsup_{r\rightarrow\infty}B(t){\bf z}^{(k_{r})}(t)
\leq {\bf c}(t)+\mbox{\boldmath $\epsilon$}+\int_{0}^{t} K(t,s){\bf z}^{(\epsilon)}(s)ds
\mbox{ a.e. in $[0,T]$.}
\end{equation}
This shows that ${\bf z}^{(\epsilon)}$ is a feasible solution of $(\mbox{CLP}_{\epsilon})$.
Let $N_{0j}=\{t\in [0,T]:z_{j}^{(\epsilon)}(t)<0\}$ and $N_{0}=\bigcup_{j=1}^{q}N_{0j}$.
Let $N_{1}$ be the subset of $[0,T]$ such that the inequality
(\ref{opteq148}) is violated. We define $N=N_{0}\cup N_{1}$.
Then from (\ref{clpeq1}) and (\ref{opteq148}), we see that the set $N$ has measure zero.
Now, we define
\begin{equation}{\label{*clpeq51}}
\bar{\bf z}^{(\epsilon)}(t)=\left\{\begin{array}{ll}
{\bf z}^{(\epsilon)}(t) & \mbox{if $t\not\in N$}\\
{\bf 0} & \mbox{if $t\in N$}.
\end{array}\right .
\end{equation}
Then, we see that $\bar{\bf z}^{(\epsilon)}(t)\geq {\bf 0}$ for all $t\in [0,T]$
and $\bar{\bf z}^{(\epsilon)}(t)={\bf z}^{(\epsilon)}(t)$ a.e. in $[0,T]$.
For $t\not\in N$, from (\ref{opteq148}), we have
\begin{equation}{\label{*clpeq52}}
B(t)\bar{\bf z}^{(\epsilon)}(t)=B(t){\bf z}^{(\epsilon)}(t)
\leq {\bf c}(t)+\mbox{\boldmath $\epsilon$}+\int_{0}^{t} K(t,s){\bf z}^{(\epsilon)}(s)ds
={\bf c}(t)+\mbox{\boldmath $\epsilon$}+\int_{0}^{t} K(t,s)\bar{\bf z}^{(\epsilon)}(s)ds.
\end{equation}
This shows that $\bar{\bf z}^{(\epsilon)}$ is a feasible solution of $(\mbox{CLP}_{\epsilon})$.

To prove part (ii), under the assumptions of ${\bf c}(t)$ and $K(t,s)$,
it is obvious that the problem $(\mbox{CLP}_{\epsilon})$ is feasible
with the trivial feasible solution ${\bf z}(t)={\bf 0}$ for all $t\in [0,T]$.
We consider $\bar{\bf z}^{(\epsilon)}$ defined in (\ref{*clpeq51}).
For $t\in N$, we have $B(t)\bar{\bf z}^{(\epsilon)}(t)={\bf 0}$.
Since ${\bf z}^{(\epsilon)}(t)\geq {\bf 0}$ a.e. in $[0,T]$ and
$K(t_{0},s)\geq {\bf 0}$ a.e. in $[0,T]$ for each fixed $t_{0}\in [0,T]$,
we obtain
\[B(t)\bar{\bf z}^{(\epsilon)}(t)={\bf 0}\leq {\bf c}(t)+\mbox{\boldmath $\epsilon$}
+\int_{0}^{t} K(t,s){\bf z}^{(\epsilon)}(s)ds
={\bf c}(t)+\mbox{\boldmath $\epsilon$}+\int_{0}^{t} K(t,s)\bar{\bf z}^{(\epsilon)}(s)ds.\]
By referring to (\ref{*clpeq52}) for $t\not\in N$,
we see that $\bar{\bf z}^{(\epsilon)}(t)$ satisfies all the constraints of
primal problem $(\mbox{CLP}_{\epsilon})$ for all $t\in [0,T]$. This completes the proof.
\end{Proof}

\begin{Thm}{\label{optt120}}
Assume that $B(t)\geq {\bf 0}$ a.e. in $[0,T]$.
For any $\epsilon\geq 0$, the following results hold.
\begin{enumerate}
\item [{\em (i)}] Suppose that the problem $(\mbox{\em CLP}_{\epsilon})$
is feasible, and that the feasible set ${\cal Z}_{\epsilon}$ of
$(\mbox{\em CLP}_{\epsilon})$ is uniformly essentially bounded.
Then, there exists an optimal solution
$\bar{\bf z}^{(\epsilon)}$ of $(\mbox{\em CLP}_{\epsilon})$ such that
$\bar{\bf z}^{(\epsilon)}(t)\geq {\bf 0}$ for all $t\in [0,T]$.

\item [{\em (ii)}] Suppose that ${\bf c}(t)\geq {\bf 0}$ for all $t\in [0,T]$,
$K(t_{0},s)\geq {\bf 0}$ a.e. in $[0,T]$ for each fixed $t_{0}\in [0,T]$,
and that the feasible set ${\cal Z}_{\epsilon}$ of $(\mbox{\em CLP}_{\epsilon})$
is uniformly essentially bounded. Then, there exists a common optimal solution
$\bar{\bf z}^{(\epsilon)}$ of $(\mbox{\em CLP}_{\epsilon})$
and $(\mbox{\em CLP}_{\epsilon}^{*})$
such that both problems have the same optimal objective values.
\end{enumerate}
\end{Thm}
\begin{Proof}
To prove part (i), we define
\[M=\sup_{{\bf z}\in {\cal Z}_{\epsilon}}\int_{0}^{T}{\bf a}^{\top}(t){\bf z}(t)dt.\]
Then, there exists a sequence $\{{\bf z}^{(k)}\}_{k=1}^{\infty}$ in ${\cal Z}_{\epsilon}$ such that
\begin{equation}{\label{opteq146}}
\lim_{k\rightarrow\infty}\int_{0}^{T}{\bf a}^{\top}(t){\bf z}^{(k)}(t)dt=M.
\end{equation}
We are going to claim that the supremum $M$ can be attained by some feasible
solution of $(\mbox{CLP}_{\epsilon})$. Since the sequence $\{{\bf z}^{(k)}\}_{k=1}^{\infty}$
is uniformly essentially bounded in $[0,T]$ by the assumption on the feasible
set ${\cal Z}_{\epsilon}$, using part (i) of Proposition~\ref{*clpr62},
there exists a subsequence $\{{\bf z}^{(k_{r})}\}_{r=1}^{\infty}$ which
weakly converges to some feasible solution ${\bf z}^{(\epsilon)}\in L_{q}^{2}[0,T]$ of $(\mbox{CLP}_{\epsilon})$, and
there exists another feasible solution $\bar{\bf z}^{(\epsilon)}$ of $(\mbox{CLP}_{\epsilon})$
such that $\bar{\bf z}^{(\epsilon)}(t)\geq {\bf 0}$ for all $t\in [0,T]$
and $\bar{\bf z}^{(\epsilon)}(t)={\bf z}^{(\epsilon)}(t)$ a.e. in $[0,T]$.
From (\ref{opteq146}), we obtain
\begin{align*}
\int_{0}^{T}{\bf a}^{\top}(t)\bar{\bf z}^{(\epsilon)}(t)dt
& =\int_{0}^{T}{\bf a}^{\top}(t){\bf z}^{(\epsilon)}(t)dt
=\sum_{j=1}^{q}\int_{0}^{T}z_{j}^{(\epsilon)}(t)a_{j}(t)dt\\
& =\lim_{r\rightarrow\infty}\sum_{j=1}^{q}\int_{0}^{T}z_{j}^{(k_{r})}(t)a_{j}(t)dt
=\lim_{r\rightarrow\infty}\int_{0}^{T}{\bf a}^{\top}(t){\bf z}^{(k_{r})}(t)dt=M.
\end{align*}
This shows that $\bar{\bf z}^{(\epsilon)}$ is an optimal solution of $(\mbox{CLP}_{\epsilon})$.

To prove part (ii), since ${\bf c}(t)\geq {\bf 0}$ for all $t\in [0,T]$ and
$K(t_{0},s)\geq {\bf 0}$ a.e. in $[0,T]$ for each fixed $t_{0}\in [0,T]$,
part (ii) of Proposition~\ref{*clpr62} says that we can take
$\bar{\bf z}^{(\epsilon)}$ as a feasible solution of $(\mbox{CLP}_{\epsilon}^{*})$.
Since the feasible set of $(\mbox{CLP}_{\epsilon}^{*})$ is contained in the
feasible set of $(\mbox{CLP}_{\epsilon})$, it follows that $\bar{\bf z}^{(\epsilon)}$
is an optimal solution of problem $(\mbox{CLP}_{\epsilon}^{*})$. This completes the proof.
\end{Proof}

Suppose that there exists a constant $\sigma >0$ such that
$\sum_{i=1}^{p}B_{ij}(t)\geq\sigma$ a.e. in $[0,T]$ for each $j=1,\cdots ,q$.
We define a real-valued function
\begin{equation}{\label{opteq153}}
\rho_{\epsilon}(t)=\frac{\tau -\epsilon}{\sigma}\cdot\exp
\left [\frac{\nu\cdot (T-t)}{\sigma}\right ]
\mbox{ for }t\in [0,T]
\end{equation}
and define $\mbox{\boldmath $\rho$}_{\epsilon}(t)$ as an $p$-dimensional
vector-valued function with all entries $\rho_{\epsilon}(t)$.
In the sequel, we are going to study the existence of optimal
solutions of $(\mbox{DCLP}_{\epsilon})$.
We first present the feasibility of dual problem.

\begin{Pro}{\label{*clpp50}}
The following statements hold true.
\begin{enumerate}
\item [{\em (i)}] Suppose that there exists a constant $\sigma >0$ such that
$\sum_{i=1}^{p}B_{ij}(t)\geq\sigma$ a.e. in $[0,T]$ for each $j=1,\cdots ,q$.
Then, the problem $(\mbox{\em DCLP}_{\epsilon})$ is feasible with the feasible
solution $\mbox{\boldmath $\rho$}_{\epsilon}$.

\item [{\em (ii)}] Suppose that there exists a constant $\sigma >0$ such that
$\sum_{i=1}^{p}B_{ij}(t)\geq\sigma$ for all $t\in [0,T]$ and for each $j=1,\cdots ,q$,
and that the function $\sum_{i=1}^{p}K_{ij}$ is bounded by $\nu$ and the function $a_{j}$
is bounded by $\tau$ for each $j=1,\cdots ,q$.
Then, the problem $(\mbox{\em DCLP}_{\epsilon}^{*})$ is feasible with the feasible
solution $\mbox{\boldmath $\rho$}_{\epsilon}$.
\end{enumerate}
\end{Pro}
\begin{Proof}
To prove part (i), from (\ref{opteq153}), we see that
\begin{equation}{\label{opteq154}}
\sigma\rho_{\epsilon}(t)=\tau -\epsilon +\nu\cdot\int_{t}^{T}\rho_{\epsilon}(s)ds
\mbox{ for $t\in [0,T]$.}
\end{equation}
For each $j=1,\cdots ,q$, using (\ref{opteq154}), we have
\begin{equation}{\label{opteq155}}
\sum_{i=1}^{p}B_{ij}(t)\rho_{\epsilon}(t)\geq\sigma\rho_{\epsilon}(t)\geq
a_{j}(t)-\epsilon +\sum_{i=1}^{p}\int_{t}^{T}
K_{ij}(s,t)\rho_{\epsilon}(s)ds\mbox{ a.e. in $[0,T]$}.
\end{equation}
This shows that $\mbox{\boldmath $\rho$}_{\epsilon}$ is a feasible solution
of $(\mbox{DCLP}_{\epsilon})$.
To prove part (ii), by applying the assumptions to (\ref{opteq155}), we obtain
\begin{equation}{\label{opteq*155}}
\sum_{i=1}^{p}B_{ij}(t)\rho_{\epsilon}(t)\geq\sigma\rho_{\epsilon}(t)\geq
a_{j}(t)-\epsilon +\sum_{i=1}^{p}\int_{t}^{T}
K_{ij}(s,t)\rho_{\epsilon}(s)ds\mbox{ for all }t\in [0,T].
\end{equation}
This completes the proof.
\end{Proof}

\begin{Lem}{\label{*clpp52}}
Let ${\bf w}^{(\epsilon)}$ be a feasible solution of problem
$(\mbox{\em DCLP}_{\epsilon})$. Then, the following statements hold true.
\begin{enumerate}
\item [{\em (i)}] There exists a feasible solution $\bar{\bf w}^{(\epsilon)}$ of
$(\mbox{\em DCLP}_{\epsilon})$
such that $\bar{\bf w}^{(\epsilon)}(t)={\bf w}^{(\epsilon)}(t)$ a.e. in $[0,T]$
and $\bar{\bf w}^{(\epsilon)}(t)\geq {\bf 0}$ for all $t\in [0,T]$.
If we further assumed that there is a vector-valued function
${\bf v}^{(\epsilon)}(t)\geq {\bf 0}$ for all $t\in [0,T]$
such that ${\bf w}^{(\epsilon)}(t)\leq
{\bf v}^{(\epsilon)}(t)$ a.e. in $[0,T]$, then
${\bf 0}\leq\bar{\bf w}^{(\epsilon)}(t)\leq{\bf v}^{(\epsilon)}(t)$ for all $t\in [0,T]$.

\item [{\em (ii)}] Suppose that there exists a constant $\sigma >0$ such that
$\sum_{i=1}^{p}B_{ij}(t)\geq\sigma$ for all $t\in [0,T]$ and for each $j=1,\cdots ,q$,
and that the function $\sum_{i=1}^{p}K_{ij}$ is bounded by $\nu$ and the function $a_{j}$
is bounded by $\tau$ for each $j=1,\cdots ,q$. If ${\bf w}^{(\epsilon)}(t)\leq
\mbox{\boldmath $\rho$}_{\epsilon}(t)$ a.e. in $[0,T]$, then there exists a feasible solution
$\bar{\bf w}^{(\epsilon)}$ of $(\mbox{\em DCLP}_{\epsilon}^{*})$
such that ${\bf 0}\leq\bar{\bf w}^{(\epsilon)}(t)\leq
\mbox{\boldmath $\rho$}_{\epsilon}(t)$ for all $t\in [0,T]$ and
$\bar{\bf w}^{(\epsilon)}(t)={\bf w}^{(\epsilon)}(t)$ a.e. in $[0,T]$.
\end{enumerate}
\end{Lem}
\begin{Proof}
To prove part (i), we begin by observing that ${\bf w}^{(\epsilon)}(t)\geq {\bf 0}$
a.e in $[0,T]$ and
\begin{equation}{\label{*clpeq54}}
B^{\top}(t){\bf w}^{(\epsilon)}(t)\geq
{\bf a}(t)-\mbox{\boldmath $\epsilon$}+\int_{t}^{T}
K^{\top}(s,t){\bf w}^{(\epsilon)}(s)ds\mbox{ a.e. in $[0,T]$}.
\end{equation}
Let $N_{0i}=\{t\in [0,T]:w_{i}^{(\epsilon)}(t)<0\}$
and $N_{0}=\bigcup_{i=1}^{p}N_{0i}$.
Let $N_{1}$ be the subset of $[0,T]$ on which the inequality
(\ref{*clpeq54}) is violated, and let $N=N_{0}\cup N_{1}$.
Then, we see that the set $N$ has measure zero. Now, we define
\[\bar{\bf w}^{(\epsilon)}(t)=\left\{\begin{array}{ll}
{\bf w}^{(\epsilon)}(t) & \mbox{if $t\not\in N$}\\
{\bf 0} & \mbox{if $t\in N$}.
\end{array}\right .\]
Then, we see that $\bar{\bf w}^{(\epsilon)}(t)\geq {\bf 0}$ for all $t\in [0,T]$ and
$\bar{\bf w}^{(\epsilon)}(t)={\bf w}^{(\epsilon)}(t)$ a.e. in $[0,T]$.
For $t\not\in N$, from (\ref{*clpeq54}), we have
\begin{align}
B^{\top}(t)\bar{\bf w}^{(\epsilon)}(t)
& =B^{\top}(t){\bf w}^{(\epsilon)}(t)\nonumber\\
& \geq {\bf a}(t)-\mbox{\boldmath $\epsilon$}
+\int_{t}^{T} K^{\top}(s,t){\bf w}^{(\epsilon)}(s)ds
={\bf a}(t)-\mbox{\boldmath $\epsilon$}
+\int_{t}^{T} K^{\top}(s,t)\bar{\bf w}^{(\epsilon)}(s)ds.\label{extdclp1}
\end{align}
This shows that $\bar{\bf w}^{(\epsilon)}$ is a feasible solution of the dual
problem $(\mbox{DCLP}_{\epsilon})$.

Now, we assume that ${\bf w}^{(\epsilon)}(t)\leq{\bf v}^{(\epsilon)}(t)$ a.e. in $[0,T]$.
Let $N_{0}$ and $N_{1}$ be the subsets of $[0,T]$ defined above,
$N_{2i}=\{t\in [0,T]:w_{i}^{(\epsilon)}(t)>v_{i}^{(\epsilon)}(t)\}$,
$N_{2}=\bigcup_{i=1}^{p}N_{2i}$, and $\hat{N}=N_{0}\cup N_{1}\cup N_{2}$.
Then, the set $\hat{N}$ has measure zero. Now, we define
\begin{equation}{\label{*clpeq55}}
\bar{\bf w}^{(\epsilon)}(t)=\left\{\begin{array}{ll}
{\bf w}^{(\epsilon)}(t) & \mbox{if $t\not\in\hat{N}$}\\
{\bf v}^{(\epsilon)}(t) & \mbox{if $t\in\hat{N}$}.
\end{array}\right .
\end{equation}
Then, we see that ${\bf 0}\leq\bar{\bf w}^{(\epsilon)}(t)\leq {\bf v}^{(\epsilon)}(t)$
for all $t\in [0,T]$ and $\bar{\bf w}^{(\epsilon)}(t)
={\bf w}^{(\epsilon)}(t)$ a.e. in $[0,T]$. For $t\not\in\hat{N}$, we have $t\not\in N$.
Using (\ref{extdclp1}), it follows that $\bar{\bf w}^{(\epsilon)}$ is a feasible
solution of the dual problem $(\mbox{DCLP}_{\epsilon})$.

To prove part (ii), from (\ref{opteq*155}), we can also obtain the following inequality
\begin{equation}{\label{clpeq266}}
B^{\top}(t)\mbox{\boldmath $\rho$}_{\epsilon}(t)\geq {\bf a}(t)
-\mbox{\boldmath $\epsilon$}+\int_{t}^{T}
K^{\top}(s,t)\mbox{\boldmath $\rho$}_{\epsilon}(s)ds\mbox{ for all $t\in [0,T]$.}
\end{equation}
We take $\bar{\bf w}^{(\epsilon)}(t)$ as
defined in (\ref{*clpeq55}) by substituting ${\bf v}^{(\epsilon)}$ for
$\mbox{\boldmath $\rho$}_{\epsilon}$. Then, we see that ${\bf 0}\leq
\bar{\bf w}^{(\epsilon)}(t)\leq\mbox{\boldmath $\rho$}_{\epsilon}(t)$
for all $t\in [0,T]$. For $t\in\hat{N}$, using (\ref{clpeq266}), we obtain
\begin{align*}
B^{\top}(t)\bar{\bf w}^{(\epsilon)}(t)
& =B^{\top}(t)\mbox{\boldmath $\rho$}_{\epsilon}(t)
\geq {\bf a}(t)-\mbox{\boldmath $\epsilon$}
+\int_{t}^{T} K^{\top}(s,t)\mbox{\boldmath $\rho$}_{\epsilon}(s)ds\\
& \geq {\bf a}(t)-\mbox{\boldmath $\epsilon$}
+\int_{t}^{T} K^{\top}(s,t)\bar{\bf w}^{(\epsilon)}(s)ds.
\end{align*}
For $t\not\in\hat{N}$, the argument of part (i) is still valid.
This shows that $\bar{\bf w}^{(\epsilon)}$ satisfies the constraints of
$(\mbox{DCLP}_{\epsilon})$ for all $t\in [0,T]$, and the proof is complete.
\end{Proof}

\begin{Lem}{\label{*clpr302}}
$K(t,s)\geq {\bf 0}$ a.e. on $[0,T]\times [0,T]$ if and only if the subset
\[N_{K}=\left\{t_{0}\in [0,T]:K(s,t_{0})\not\geq {\bf 0}
\mbox{ a.e. on $[0,T]$}\right\}\]
has measure zero; that is, for each fixed $t_{0}\in [0,T]\setminus N_{K}$,
$K(s,t_{0})\geq {\bf 0}$ a.e. on $[0,T]$.
\end{Lem}
\begin{Proof}
Suppose that $K(t,s)\geq {\bf 0}$ a.e. on $[0,T]\times [0,T]$.
We are going to prove it by contradiction. Assume that $\mu (N_{K})\neq 0$.
For each fixed $t_{0}\in N_{K}$, the following set
\[\{s\in [0,T]:K(s,t_{0})\geq {\bf 0}\}\]
has measure zero, which also says that the following set
\[M_{t_{0}}\equiv\{s\in [0,T]:K(s,t_{0})\not\geq {\bf 0}\}\]
is not measure zero. Let
\[M=\bigcup_{t_{0}\in N_{K}}M_{t_{0}}.\]
Then, we have $\mu (M)\neq 0$. For each $(t,s)\in M\times N_{K}$,
we see that $K(t,s)\not\geq {\bf 0}$.
Since $(\mu\times\mu )(M\times N_{K})=\mu (M)\cdot\mu (N_{K})\neq0$,
this contradicts $K(t,s)\geq {\bf 0}$ a.e. on $[0,T]\times [0,T]$.

For the converse, let
\[{\cal N}=\left\{(s,t)\in [0,T]\times [0,T]:K(s,t)\not\geq {\bf 0}\right\}.\]
Assume that $(\mu\times\mu )({\cal N})>0$.
We are going to lead to a contradiction. It is well-know that the
Lebesgue measure $(\mu\times\mu )({\cal N})$ is equal to the inner measure given by
\[0<(\mu\times\mu )({\cal N}) =\sup_{(\bigcup_{k}{\cal R}_{k})\subseteq {\cal N}}
\sum_{k}m({\cal R}_{k}),\]
where the union is a countable union, each ${\cal R}_{k}$ is a
rectangle of $[0,T]\times [0,T]$ and
$m({\cal R}_{k})$ is the area of the rectangle ${\cal R}_{k}$.
Of course, we have $m({\cal R}_{k})=(\mu\times\mu )({\cal R}_{k})$.
In this case, there exists a rectangle ${\cal R}_{k_{0}}\subseteq
{\cal N}$ such that $0<m({\cal R}_{k_{0}})\leq(\mu\times\mu )({\cal N})$.
Suppose that ${\cal R}_{k_{0}}=R_{k_{0}}^{(1)}\times R_{k_{0}}^{(2)}$,
where $R_{k_{0}}^{(1)}$ and $R_{k_{0}}^{(2)}$ are intervals in $[0,T]$ such that
$\mu (R_{k_{0}}^{(1)})\neq 0$ and $\mu (R_{k_{0}}^{(2)})\neq 0$.
Since $\mu (N_{K})=0$ and $\mu (R_{k_{0}}^{(2)})\neq 0$,
there exists $t_{0}^{*}\in R_{k_{0}}^{(2)}$ and $t_{0}^{*}\not\in N_{K}$ such that
$R_{k_{0}}^{(1)}\times\{t_{0}^{*}\}\subseteq {\cal N}$. This shows that
$K(s,t_{0}^{*})\not\geq {\bf 0}$ for $(s,t_{0}^{*})\in
R_{k_{0}}^{(1)}\times\{t_{0}^{*}\}$, which contradicts
$K(s,t_{0}^{*})\geq {\bf 0}$ a.e. on $[0,T]$, since $t_{0}^{*}\not\in N_{K}$
and $\mu (R_{k_{0}}^{(1)})\neq 0$. This completes the proof.
\end{Proof}

\begin{Lem}{\label{optl156}}
Suppose that the following conditions are satisfied:
\begin{itemize}
\item $K(s,t)\geq {\bf 0}$ a.e. in $[0,T]\times [0,T]$;

\item $\sum_{i=1}^{p}B_{ij}(t)>0$ a.e. in $[0,T]$ for each $j=1,\cdots ,q$;

\item there exists a constant $\sigma >0$ such that, for each $i=1,\cdots ,p$ and
$j=1,\cdots ,q$, the following statement holds true a.e. in $[0,T]$:
\begin{equation}{\label{dclp117}}
B_{ij}(t)\neq 0\mbox{ implies }B_{ij}(t)\geq\sigma .
\end{equation}
\end{itemize}
Consider the vector-valued function $\mbox{\boldmath $\rho$}_{\epsilon}$
defined in $(\ref{opteq153})$, and let ${\bf w}^{(\epsilon)}$ be a feasible solution
of problem $(\mbox{\em DCLP}_{\epsilon})$.
Then, there exist a feasible solution $\widehat{\bf w}^{(\epsilon )}$ of
$(\mbox{\em DCLP}_{\epsilon})$ such that
\begin{equation}{\label{dclp114}}
\widehat{\bf w}^{(\epsilon )}(t)\geq {\bf 0}\mbox{ a.e. in $[0,T]$}
\end{equation}
and
\begin{equation}{\label{dclp115}}
\widehat{\bf w}^{(\epsilon )}(t)\leq {\bf w}^{(\epsilon)}(t)\mbox{ and }
\widehat{\bf w}^{(\epsilon )}(t)\leq\mbox{\boldmath $\rho$}_{\epsilon}(t)
\mbox{ for all $t\in [0,T]$.}
\end{equation}
Moreover, if ${\bf w}^{(\epsilon)}$ is an optimal solution of
$(\mbox{\em DCLP}_{\epsilon})$, then $\widehat{\bf w}^{(\epsilon )}$ is
also an optimal solution of $(\mbox{\em DCLP}_{\epsilon})$.
\end{Lem}
\begin{Proof}
Under the assumption of $B(t)$, it is easy to see that $B(t)\geq {\bf 0}$ a.e.
in $[0,T]$ and $\sum_{i=1}^{p}B_{ij}(t)\geq\sigma$ a.e. in $[0,T]$
for each $j=1,\cdots ,q$. Therefore, the dual problem $(\mbox{DCLP}_{\epsilon})$
is feasible by Proposition~\ref{*clpp50}.
Since ${\bf w}^{(\epsilon)}$ is a feasible solution of
dual problem $(\mbox{DCLP}_{\epsilon})$, for each $j=1,\cdots ,q$, we have
\begin{equation}{\label{opteq150}}
\sum_{i}B_{ij}(t)w_{i}^{(\epsilon)}(t)\geq a_{j}(t)-\epsilon+\sum_{i}\int_{t}^{T}
K_{ij}(s,t)w_{i}^{(\epsilon)}(s)ds\mbox{ a.e. in $[0,T]$}.
\end{equation}
Now, for $t\in [0,T]$, we define
\[\widehat{w}_{i}^{(\epsilon )}(t)=\min\left\{w_{i}^{(\epsilon)}(t),
\rho_{\epsilon}(t)\right\}.\]
It is obvious that (\ref{dclp114}) and (\ref{dclp115}) are satisfied.
On the other hand, from (\ref{opteq150}) we also obtain
\begin{equation}{\label{opteq151}}
\sum_{i}B_{ij}(t)w_{i}^{(\epsilon)}(t)\geq a_{j}(t)-\epsilon+\sum_{i}\int_{t}^{T}
K_{ij}(s,t)\widehat{w}_{i}^{(\epsilon )}(s)ds\mbox{ a.e. in $[0,T]$}.
\end{equation}

Let $\bar{N}_{1i}=\{t\in [0,T]:\widehat{w}_{i}^{(\epsilon )}(t)<0\}$ and
$\bar{N}_{1}=\bigcup_{i=1}^{p}\bar{N}_{1i}$.
Let $\bar{N}_{2ij}=\{t\in [0,T]:B_{ij}(t)<0\}$ and
$\bar{N}_{2}=\bigcup_{i=1}^{p}\bigcup_{j=1}^{q}\bar{N}_{2ij}$.
Let $\bar{N}_{3}$ be the subset of $[0,T]$ on which the statement (\ref{dclp117})
is violated. Let
\[\bar{N}=\bar{N}_{1}\cup\bar{N}_{2}\cup\bar{N}_{3}\cup N_{K},\]
where $N_{K}$ is defined in Lemma~\ref{*clpr302}. Then $\bar{N}$ has measure zero.
Let $N_{0}$ and $N_{1}$ be the subsets of $[0,T]$ on which
the inequalities (\ref{opteq155}) and (\ref{opteq151}) are violated, respectively. We take
\[N=N_{0}\cup N_{1}\cup\bar{N}.\]
Then, the set $N$ has measure zero. For any fixed $t\in [0,T]\setminus N$,
we define the index sets $I_{\leq}=\{i:w_{i}^{(\epsilon)}(t)\leq\rho_{\epsilon}(t)\}$
and $I_{>}=\{i:w_{i}^{(\epsilon)}(t)>\rho_{\epsilon}(t)\}$, and consider
\[\sum_{i}B_{ij}(t)\widehat{w}_{i}^{(\epsilon )}(t)=\sum_{i\in I_{\leq}}
B_{ij}(t)\widehat{w}_{i}^{(\epsilon )}(t)
+\sum_{i\in I_{>}}B_{ij}(t)\widehat{w}_{i}^{(\epsilon )}(t).\]
Then, we have the following three cases.
\begin{itemize}
\item Suppose that $I_{>}=\emptyset$ (i.e., the second sum is zero).
Then, we see that $w_{i}^{(\epsilon)}(t)=\widehat{w}_{i}^{(\epsilon )}(t)$ for all $i$.
Therefore, from (\ref{opteq151}) we have
\[\sum_{i}B_{ij}(t)\widehat{w}_{i}^{(\epsilon )}(t)=\sum_{i}B_{ij}(t)
w_{i}^{(\epsilon)}(t)\geq a_{j}(t)-\epsilon+\sum_{i}\int_{t}^{T}
K_{ij}(s,t)\widehat{w}_{i}^{(\epsilon )}(s)ds.\]

\item Suppose that $I_{>}\neq\emptyset$ and $B_{ij}(t)=0$ for all $i\in I_{>}$.
Then, by (\ref{opteq151}), we also have
\[\sum_{i}B_{ij}(t)\widehat{w}_{i}^{(\epsilon )}(t)
=\sum_{i}B_{ij}(t)w_{i}^{(\epsilon)}(t)
\geq a_{j}(t)-\epsilon+\sum_{i}\int_{t}^{T}
K_{ij}(s,t)\widehat{w}_{i}^{(\epsilon )}(s)ds.\]

\item Suppose that $I_{>}\neq\emptyset$, and that there exists $i_{0}\in I_{>}$ with
$B_{i_{0}j}(t)\neq 0$, i.e., $B_{i_{0}j}(t)\geq\sigma$ by the
assumption on $B(t)$ (since $t\not\in\bar{N}_{3}$).
Since $t\not\in\bar{N}_{0}$ and $t\not\in\bar{N}_{2}$, it follows that
$\widehat{w}_{i}^{(\epsilon )}(t)\geq 0$ and $B_{ij}(t)\geq 0$
for $i=1,\cdots ,p$ and $j=1,\cdots ,q$. Therefore, we have
\begin{equation}{\label{opteq156}}
\sum_{i}B_{ij}(t)\widehat{w}_{i}^{(\epsilon )}(t)\geq
\sum_{i\in I_{>}}B_{ij}(t)\widehat{w}_{i}^{(\epsilon )}(t)=
\sum_{i\in I_{>}}B_{ij}(t)\rho_{\epsilon}(t)
\geq B_{i_{0}j}(t)\rho_{\epsilon}(t)\geq\sigma\rho_{\epsilon}(t).
\end{equation}
Since the fixed $t\not\in N_{K}$, it follows that $K_{ij}(s,t)\geq 0$ a.e. in
$[0,T]$ by Lemma~\ref{*clpr302}, Since $t\not\in N_{0}$,
using (\ref{opteq155}) and (\ref{opteq156}), we obtain
\begin{align*}
\sum_{i}B_{ij}(t)\widehat{w}_{i}^{(\epsilon )}(t)
& \geq a_{j}(t)-\epsilon+\sum_{i}\int_{t}^{T}
K_{ij}(s,t)\rho_{\epsilon}(s)ds\\
& \geq a_{j}(t)-\epsilon+\sum_{i}\int_{t}^{T}
K_{ij}(s,t)\widehat{w}_{i}^{(\epsilon )}(s)ds.
\end{align*}
\end{itemize}
Therefore, we conclude that
\[\sum_{i}B_{ij}(t)\widehat{w}_{i}^{(\epsilon )}(t)\geq a_{j}(t)-\epsilon
+\sum_{i}\int_{t}^{T}
K_{ij}(s,t)\widehat{w}_{i}^{(\epsilon )}(s)ds\mbox{ a.e in $[0,T]$.}\]
This shows that $\widehat{\bf w}^{(\epsilon )}$ is a feasible solution of
$(\mbox{DCLP}_{\epsilon})$.
Suppose that ${\bf w}^{(\epsilon)}$ is an optimal solution of
$(\mbox{DCLP}_{\epsilon})$. Since $(\mbox{DCLP}_{\epsilon})$
is a minimization problem and $\widehat{\bf w}^{(\epsilon )}(t)\leq
{\bf w}^{(\epsilon)}(t)$ for all $t\in [0,T]$, we have
\[\int_{0}^{T}\left ({\bf c}^{(\epsilon)}(t)\right )^{\top}\widehat{\bf w}^{(\epsilon )}(t)dt
\leq\int_{0}^{T}\left ({\bf c}^{(\epsilon)}(t)\right )^{\top}{\bf w}^{(\epsilon)}(t)dt
\leq\int_{0}^{T}\left ({\bf c}^{(\epsilon)}(t)\right )^{\top}\widehat{\bf w}^{(\epsilon )}(t)dt,\]
which says that $\widehat{\bf w}^{(\epsilon )}$ is an optimal solution of
$(\mbox{DCLP}_{\epsilon})$. This completes the proof.
\end{Proof}

\begin{Rem}{\label{*clpr97}}
{\em
We see that if the assumption regarding the time-dependent matrix $B(t)$ in
Lemma~\ref{optl156} is satisfied, then $\sum_{i=1}^{p}B_{ij}(t)\geq\sigma$
a.e. in $[0,T]$ for each $j=1,\cdots ,q$, which says that
the assumption of Proposition~\ref{optt120*} regarding
the time-dependent matrix $B(t)$ is also satisfied by taking
$\lambda_{i}(t)=1$ for all $i=1,\cdots ,p$ and $t\in [0,T]$.
In other words, the conclusions of Proposition~\ref{optt120*} are available
when the assumption in Lemma~\ref{optl156} is satisfied.
}\end{Rem}

\begin{Pro}{\label{*clpp48*}}
Suppose that the following conditions are satisfied:
\begin{itemize}
\item $K(s,t)\geq {\bf 0}$ a.e. in $[0,T]\times [0,T]$;

\item $\sum_{i=1}^{p}B_{ij}(t)>0$ a.e. in $[0,T]$ for each $j=1,\cdots ,q$;

\item there exists a constant $\sigma >0$ such that, for each $i=1,\cdots ,p$ and
$j=1,\cdots ,q$, the following statement holds true a.e. in $[0,T]$:
\[B_{ij}(t)\neq 0\mbox{ implies }B_{ij}(t)\geq\sigma .\]
\end{itemize}
Consider the sequence $\{\epsilon_{k}\}_{k=1}^{\infty}$ with
$\epsilon_{k}\rightarrow 0+$ as $k\rightarrow\infty$. For each $\epsilon_{k}$,
let ${\bf w}^{(\epsilon_{k})}$ be a feasible solution of problem
$(\mbox{\em DCLP}_{\epsilon_{k}})$. Then, for each $\epsilon_{k}$, there exists a
feasible solution $\widehat{\bf w}^{(\epsilon_{k})}$
of problem $(\mbox{\em DCLP}_{\epsilon_{k}})$ such that
the following properties hold true.
\begin{enumerate}
\item [{\em (i)}] The sequence $\{\widehat{\bf w}^{(\epsilon_{k})}\}_{k=1}^{\infty}$
is uniformly bounded.

\item [{\em (ii)}] $\widehat{\bf w}^{(\epsilon_{k})}(t)\leq
{\bf w}^{(\epsilon_{k})}(t)$ for all $t\in [0,T]$ and, for each $i=1,\cdots ,p$,
\[\widehat{w}_{i}^{(\epsilon_{k})}(t)\geq 0\mbox{ a.e. in }[0,T]\]
and
\[\widehat{w}_{i}^{(\epsilon_{k})}(t)\leq\frac{\tau -\epsilon_{k}}{\sigma}\cdot\exp
\left [\frac{\nu\cdot (T-t)}{\sigma}\right ]\leq
\frac{\tau}{\sigma}\cdot\exp\left (\frac{\nu\cdot T}{\sigma}\right )
\mbox{ for all $t\in [0,T]$}.\]

\item [{\em (iii)}] There exists a subsequence
$\{\widehat{\bf w}^{(\epsilon_{k_{r}})}\}_{r=1}^{\infty}$
which weakly converges to some feasible solution ${\bf w}^{(0)}\in L_{p}^{2}[0,T]$ of problem
$(\mbox{\em DCLP}_{0})=(\mbox{\em DCLP})$. Moreover, there is also another
feasible solution $\bar{\bf w}$ of problem {\em (DCLP)} such that
$\bar{\bf w}(t)={\bf w}^{(0)}(t)$ a.e. in $[0,T]$ and, for each $i=1,\cdots ,p$,
\begin{equation}{\label{*clpeq65}}
0\leq\bar{w}_{i}(t)\leq\frac{\tau}{\sigma}\cdot\exp
\left [\frac{\nu\cdot (T-t)}{\sigma}\right ]\leq\frac{\tau}{\sigma}\cdot\exp
\left (\frac{\nu\cdot T}{\sigma}\right )
\mbox{ for all $t\in [0,T]$}.
\end{equation}
We further assume that the conditions regarding the time-dependent matrix $B(t)$
are satisfied for all $t\in [0,T]$,
and that the function $\sum_{i=1}^{p}K_{ij}$ is bounded by $\nu$ and the function $a_{j}$
is bounded by $\tau$ for each $j=1,\cdots ,q$.
Then $\bar{\bf w}(t)$ can be taken as a feasible solution of $(\mbox{\em DCLP}^{*})$.
\end{enumerate}
\end{Pro}
\begin{Proof}
By Lemma~\ref{optl156}, there exists a sequence
$\{\widehat{\bf w}^{(\epsilon_{k})}\}_{k=1}^{\infty}$ of feasible solutions
of problems $(\mbox{DCLP}_{\epsilon_{k}})$ such that
$\widehat{\bf w}^{(\epsilon_{k})}(t)\leq {\bf w}^{(\epsilon_{k})}(t)$
for all $t\in [0,T]$ and, for each $i=1,\cdots ,p$,
\begin{equation}{\label{dclp113}}
\widehat{w}_{i}^{(\epsilon_{k})}(t)\geq 0\mbox{ a.e. in }[0,T]
\end{equation}
and
\begin{equation}{\label{*clpeq64}}
\widehat{w}_{i}^{(\epsilon_{k})}(t)\leq\rho_{\epsilon_{k}}(t)
=\frac{\tau -\epsilon_{k}}{\sigma}\cdot\exp
\left [\frac{\nu\cdot (T-t)}{\sigma}\right ]\leq\frac{\tau}{\sigma}\cdot\exp
\left (\frac{\nu\cdot T}{\sigma}\right )\mbox{ for all $t\in [0,T]$,}
\end{equation}
which says that the sequence $\{\widehat{\bf w}^{(\epsilon_{k})}\}_{k=1}^{\infty}$
is uniformly bounded in $[0,T]$. This proves parts (i) and (ii).
Now, using Lemma~\ref{optl138}, there exists a subsequence
$\{\widehat{\bf w}^{(\epsilon_{k_{r}})}\}_{r=1}^{\infty}$ which weakly converges
to some ${\bf w}^{(0)}\in L_{p}^{2}[0,T]$. Using Lemma~\ref{optl157}, we have
\begin{equation}{\label{opteq223}}
{\bf 0}\leq\liminf_{r\rightarrow\infty}
\widehat{\bf w}^{(\epsilon_{k_{r}})}(t)\leq {\bf w}^{(0)}(t)
\leq\limsup_{r\rightarrow\infty}
\widehat{\bf w}^{(\epsilon_{k_{r}})}(t)\mbox{ a.e. in $[0,T]$}.
\end{equation}
Since $\{\widehat{\bf w}^{(\epsilon_{k_{r}})}\}_{r=1}^{\infty}$ are feasible
solutions of problems $(\mbox{DCLP}_{\epsilon_{k_{r}}})$, we have
\begin{equation}{\label{*clpeq67}}
B^{\top}(t) \widehat{\bf w}^{(\epsilon_{k_{r}})}(t)\geq {\bf a}(t)
-\mbox{\boldmath $\epsilon$}_{k_{r}}+\int_{t}^{T}
K^{\top}(s,t)\widehat{\bf w}^{(\epsilon_{k_{r}})}(s)ds\mbox{ a.e. in $[0,T]$}.
\end{equation}
By taking the limit inferior and using the weak convergence from (\ref{*clpeq67}),
since $B(t)\geq {\bf 0}$ a.e. in $[0,T]$, using (\ref{opteq223}), we obtain
\[B^{\top}(t){\bf w}^{(0)}(t)\geq\liminf_{r\rightarrow\infty}
B^{\top}(t)\widehat{\bf w}^{(\epsilon_{k_{r}})}(t)
\geq {\bf a}(t)+\int_{t}^{T} K^{\top}(s,t){\bf w}^{(0)}(s)ds\mbox{ a.e. in $[0,T]$}.\]
This shows that ${\bf w}^{(0)}$ is a feasible solution of problem (DCLP).
Using (\ref{opteq223}) and (\ref{*clpeq64}), we have
\[0\leq w_{i}^{(0)}(t)\leq\frac{\tau}{\sigma}\cdot\exp
\left [\frac{\nu\cdot (T-t)}{\sigma}\right ]=\rho_{0}(t)\mbox{ a.e. in $[0,T]$
for each $i=1,\cdots ,p$}.\]
Using part (i) of Lemma~\ref{*clpp52} by taking $\epsilon =0$ with
$\bar{w}_{i}\equiv\bar{w}_{i}^{(0)}$ and
\[v_{i}^{(0)}(t)\equiv\frac{\tau}{\sigma}\cdot\exp
\left [\frac{\nu\cdot (T-t)}{\sigma}\right ]=\rho_{0}(t),\]
we obtain the desired result.

Finally, we further assume that the conditions regarding the time-dependent matrix $B(t)$
are satisfied for all $t\in [0,T]$,
and that the function $\sum_{i=1}^{p}K_{ij}$ is bounded by $\nu$ and the function $a_{j}$
is bounded by $\tau$ for each $j=1,\cdots ,q$. Then, the desired result
follows from part (ii) of Lemma~\ref{*clpp52} by taking $\epsilon =0$. This completes the proof.
\end{Proof}

\begin{Thm}{\label{p63}}
Suppose that the following conditions are satisfied:
\begin{itemize}
\item $K(s,t)\geq {\bf 0}$ a.e. in $[0,T]\times [0,T]$;

\item $\sum_{i=1}^{p}B_{ij}(t)>0$ a.e. in $[0,T]$ for each $j=1,\cdots ,q$;

\item there exists a constant $\sigma >0$ such that, for each $i=1,\cdots ,p$ and
$j=1,\cdots ,q$, the following statement holds true a.e. in $[0,T]$:
\[B_{ij}(t)\neq 0\mbox{ implies }B_{ij}(t)\geq\sigma .\]
\end{itemize}
Then, the following results hold.
\begin{enumerate}
\item [{\em (i)}] The problem $(\mbox{\em DCLP}_{\epsilon})$ has an optimal solution
$\bar{\bf w}^{(\epsilon )}$ such that, for each $i=1,\cdots ,p$,
\begin{equation}{\label{*clpeq68}}
0\leq\bar{w}_{i}^{(\epsilon )}(t)\leq\frac{\tau}{\sigma}\cdot\exp
\left [\frac{\nu\cdot (T-t)}{\sigma}\right ]\leq\frac{\tau}{\sigma}\cdot\exp
\left (\frac{\nu\cdot T}{\sigma}\right )
\mbox{ for all $t\in [0,T]$}.
\end{equation}

\item [{\em (ii)}] We further assume that the conditions regarding the time-dependent matrix
$B(t)$ are satisfied for all $t\in [0,T]$,
and that the function $\sum_{i=1}^{p}K_{ij}$ is bounded by $\nu$ and the function $a_{j}$
is bounded by $\tau$ for each $j=1,\cdots ,q$.
Then, there exists a common optimal solution $\bar{\bf w}^{(\epsilon )}$ of problems
$(\mbox{\em DCLP}_{\epsilon})$ and $(\mbox{\em DCLP}_{\epsilon}^{*})$
such that the inequalities $(\ref{*clpeq68})$ are satisfied and
both problems have the same optimal objective values.
\end{enumerate}
\end{Thm}
\begin{Proof}
To prove part (i), using Proposition~\ref{*clpp50},
we see that problem $(\mbox{DCLP}_{\epsilon})$ is feasible, i.e., the feasible set
${\cal W}_{\epsilon}$ of problem $(\mbox{DCLP}_{\epsilon})$ is nonempty. Therefore, if we define
\[M=\inf_{{\bf w}^{(\epsilon)}\in {\cal W}_{\epsilon}}\int_{0}^{T}
{\bf c}^{\top}(t){\bf w}^{(\epsilon)}(t)dt.\]
Then, there exists a sequence $\{{\bf w}^{(k)}\}_{k=1}^{\infty}$ in ${\cal W}_{\epsilon}$
such that
\begin{equation}{\label{opteq146*}}
\lim_{k\rightarrow\infty}\int_{0}^{T}{\bf c}^{\top}(t){\bf w}^{(k)}(t)dt=M.
\end{equation}
By Lemma~\ref{optl156}, there exists a sequence $\{\widehat{\bf
w}^{(k)}\}_{k=1}^{\infty}$ of feasible solutions of problems
$(\mbox{DCLP}_{\epsilon})$ such that $\widehat{\bf w}^{(k)}(t)\leq
{\bf w}^{(k)}(t)$ for all $t\in [0,T]$ and, for each $i=1,\cdots ,p$,
\begin{equation}{\label{dclp113}}
\widehat{w}_{i}^{(k)}(t)\geq 0\mbox{ a.e. in }[0,T]
\end{equation}
and
\begin{equation}{\label{*2clpeq64}}
\widehat{w}_{i}^{(k)}(t)\leq\rho_{\epsilon}(t)
=\frac{\tau -\epsilon}{\sigma}\cdot\exp
\left [\frac{\nu\cdot (T-t)}{\sigma}\right ]\leq\frac{\tau}{\sigma}\cdot\exp
\left (\frac{\nu\cdot T}{\sigma}\right )\mbox{ for all $t\in [0,T]$,}
\end{equation}
which says that the sequence $\{\widehat{\bf
w}^{(k)}\}_{k=1}^{\infty}$ is uniformly bounded in $[0,T]$. Now,
using Lemma~\ref{optl138}, there exists a subsequence
$\{\widehat{\bf w}^{(k_{r})}\}_{r=1}^{\infty}$ which weakly
converges to some ${\bf w}^{(\epsilon)}\in L_{p}^{2}[0,T]$. Using Lemma~\ref{optl157}, we have
\begin{equation}{\label{2opteq223}}
{\bf 0}\leq\liminf_{r\rightarrow\infty}
\widehat{\bf w}^{(k_{r})}(t)\leq {\bf w}^{(\epsilon)}(t)
\leq\limsup_{r\rightarrow\infty}
\widehat{\bf w}^{(k_{r})}(t)\mbox{ a.e. in $[0,T]$}.
\end{equation}
Since $\{\widehat{\bf w}^{(k_{r})}\}_{r=1}^{\infty}$ are feasible
solutions of problems $(\mbox{DCLP}_{\epsilon})$, we have
\begin{equation}{\label{*2clpeq67}}
B^{\top}(t) \widehat{\bf w}^{(k_{r})}(t)\geq {\bf a}(t)
-\mbox{\boldmath $\epsilon$}+\int_{t}^{T}
K^{\top}(s,t)\widehat{\bf w}^{(k_{r})}(s)ds\mbox{ a.e. in $[0,T]$}.
\end{equation}
By taking the limit inferior and using the weak convergence from
(\ref{*2clpeq67}), since $B(t)\geq {\bf 0}$ a.e. in $[0,T]$, using (\ref{2opteq223}), we obtain
\[B^{\top}(t){\bf w}^{(\epsilon)}(t)\geq\liminf_{r\rightarrow\infty}
B^{\top}(t)\widehat{\bf w}^{(k_{r})}(t)\geq {\bf a}(t)-\mbox{\boldmath $\epsilon$}
+\int_{t}^{T} K^{\top}(s,t){\bf w}^{(\epsilon)}(s)ds\mbox{ a.e. in $[0,T]$}.\]
This shows that ${\bf w}^{(\epsilon)}$ is a feasible solution of
problem $(\mbox{DCLP}_{\epsilon})$. Using (\ref{2opteq223}) and (\ref{*2clpeq64}), we have
\[0\leq w_{i}^{(\epsilon)}(t)\leq\frac{\tau}{\sigma}\cdot\exp
\left [\frac{\nu\cdot (T-t)}{\sigma}\right ]\mbox{ a.e. in $[0,T]$ for each $i=1,\cdots ,p$}.\]
Using part (i) of Lemma~\ref{*clpp52} by taking
\[v_{i}^{(\epsilon)}(t)\equiv\frac{\tau}{\sigma}\cdot\exp
\left [\frac{\nu\cdot (T-t)}{\sigma}\right ],\]
we obtain $\bar{\bf w}^{(\epsilon)}(t)={\bf w}^{(\epsilon)}(t)$ a.e. in $[0,T]$
for some feasible solution $\bar{\bf w}^{(\epsilon)}$ of problem $(\mbox{DCLP}_{\epsilon})$
satisfying (\ref{*clpeq68}). Therefore, using the weak convergence, we have
\begin{align*}
\int_{0}^{T}{\bf c}^{\top}(t)\bar{\bf w}^{(\epsilon)}(t)dt
& =\int_{0}^{T}{\bf c}^{\top}(t){\bf w}^{(\epsilon)}(t)dt
=\lim_{r\rightarrow\infty}\int_{0}^{T}
{\bf c}^{\top}(t)\widehat{\bf w}^{(k_{r})}(t)dt\\
& \leq\lim_{r\rightarrow\infty}\int_{0}^{T}
{\bf c}^{\top}(t){\bf w}^{(k_{r})}(t)dt=M.
\end{align*}
Since $\{\widehat{\bf w}^{(k_{r})}\}_{r=1}^{\infty}$ are feasible solutions of the
minimization problem $(\mbox{DCLP}_{\epsilon})$, we have
\[\int_{0}^{T}{\bf c}^{\top}(t)\bar{\bf w}^{(\epsilon)}(t)dt
=\int_{0}^{T}{\bf c}^{\top}(t){\bf w}^{(\epsilon)}(t)dt
=\lim_{r\rightarrow\infty}\int_{0}^{T}
{\bf c}^{\top}(t)\widehat{\bf w}^{(k_{r})}(t)dt\geq M.\]
This shows that $\bar{\bf w}^{(\epsilon)}$ is an optimal solution of
problem $(\mbox{DCLP}_{\epsilon})$ such that the inequalities (\ref{*clpeq68}) are satisfied.

To prove part (ii), under the further assumptions, using part (ii) of Lemma~\ref{*clpp52},
we can see that $\bar{\bf w}^{(\epsilon)}$ is also a feasible solution of
$(\mbox{DCLP}_{\epsilon}^{*})$.
Since the feasible set of $(\mbox{DCLP}_{\epsilon}^{*})$ is contained in the
feasible set of $(\mbox{DCLP}_{\epsilon})$, we conclude that $\bar{\bf w}^{(\epsilon)}$
is also an optimal solution of $(\mbox{DCLP}_{\epsilon}^{*})$. This completes the proof.
\end{Proof}

\section{Discretized Problems}

Now, we are going to consider the discretized versions of problems
(CLP) and (DCLP). Let us consider the Lebesgue measure $\mu$ on
$[0,T]$. Then, we have $\mu ([0,T])=T$. According to the constraints of problems
(CLP) and (DCLP), we see that there is a subset ${\cal T}$ of $[0,T]$ such that
all the constraints of (CLP) and (DCLP) are satisfied
for all $t\in {\cal T}$ and $\mu ({\cal T})=T$. Let
\[{\cal P}=\left\{0=t_{0},t_{1},t_{2},\cdots ,t_{N}=T\right\}\]
be a partition of $[0,T]$ such that $t_{u}\in {\cal T}$ for all
$u=1,\cdots ,N-1$. In this case, all the constraints of
(CLP) and (DCLP) are satisfied for all $t\in {\cal P}$.

\begin{Rem}{\label{clpr42}}
{\em
Suppose that some conditions are satisfied a.e. in $[0,T]$.
We can also construct a subset $\hat{\cal T}$ of $[0,T]$ such that
$\mu (\hat{\cal T})=T$ and define a new partition $\hat{\cal P}$ of $[0,T]$
such that these conditions are satisfied for all $t\in\hat{\cal P}$.
For example, some of these conditions are listed below:
\begin{itemize}
\item the essential boundedness shown in (\ref{clpeq8})-(\ref{clpeq3});

\item ${\bf c}(t)\geq {\bf 0}$ a.e in $[0,T]$;

\item $K(t,s)\geq {\bf 0}$ a.e. on $[0,T]\times [0,T]$;

\item there exists a constant $\sigma >0$ such that, for each $i=1,\cdots ,p$ and
$j=1,\cdots ,q$, the following statement holds true a.e. in $[0,T]$:
\[B_{ij}(t)\neq 0\mbox{ implies }B_{ij}(t)\geq\sigma ;\]

\item there exist real-valued functions $\lambda_{i}$ satisfying
$0\leq\lambda_{i}(t)\leq 1$ a.e. in $[0,T]$ for $i=1,\cdots ,p$ and a constant
$\sigma >0$ satisfying
\[\min_{j=1,\cdots ,q}\left\{\sum_{i=1}^{p}\lambda_{i}(t)B_{ij}(t)
\right\}\geq\sigma\mbox{ a.e. in $[0,T]$}.\]
\end{itemize}
}\end{Rem}

Recall that $f$ is piecewise continuous on $[0,T]$ means that $f$ is continuous on
$[0,T]$ except for a finite subset of $[0,T]$.
Let $H$ be a real-valued function define on $[0,T]\times [0,T]$.
We consider the piecewise continuity for $H$ in the following sense:
there is a partition ${\cal P}_{H}=\left\{0=t_{0},t_{1},t_{2},\cdots ,t_{M}=T\right\}$
on $[0,T]$ such that the following conditions are satisfied:
\begin{itemize}
\item $H(s,t)$ is continuous on the open rectangles
$(t_{u-1},t_{u})\times (t_{v-1},t_{v})$ for $u=1,\cdots ,M$ and $v=1,\cdots ,M$;

\item for any fixed $s\not\in {\cal P}_{H}$, the single-variable function $H(s,\cdot )$ is
continuous on the open intervals $(t_{v-1},t_{v})$ for $v=1,\cdots ,M$;

\item for any fixed $t\not\in {\cal P}_{H}$, the single-variable function $H(\cdot ,t)$ is
continuous on the open intervals $(t_{u-1},t_{u})$ for $u=1,\cdots ,M$.
\end{itemize}
In order to prove the strong duality theorem, we assume further that each entry of
${\bf a}$, ${\bf c}$, $B$ and $K$ is piecewise continuous on $[0,T]$ and $[0,T]\times [0,T]$,
respectively; that is, $a_{j}$, $c_{i}$, $B_{ij}$ and $K_{ij}$ are piecewise
continuous on $[0,T]$ and $[0,T]\times [0,T]$, respectively,
for each $i=1,\cdots ,p$ and $j=1,\cdots ,q$.

Under the above assumptions, we can take the partition ${\cal P}$
such that all the discontinuities of $a_{j}$, $c_{i}$, $B_{ij}$ and $K_{ij}$
are contained in ${\cal P}$. In this case, each $a_{j}$, $c_{i}$, $B_{ij}$ and $K_{ij}$
is continuous on the open subintervals $(t_{u-1},t_{u})$ and open rectangles
$(t_{u-1},t_{u})\times (t_{v-1},t_{v})$ for $u=1,\cdots ,N$ and $v=1,\cdots ,N$, respectively.

Given a partition ${\cal P}=\{t_{0},t_{1},\cdots ,t_{N}\}$ of $[0,T]$, let
\[\parallel {\cal P}\parallel =\max_{u=1,\cdots ,N}\left (t_{u}-t_{u-1}\right )
\mbox{ satisfying }\lim_{N\rightarrow\infty}\parallel {\cal P}\parallel=0.\]
According to the above construction, we assume that the partition ${\cal P}$
satisfies the following conditions.
\begin{itemize}
\item The value $\parallel {\cal P}\parallel$ can be sufficiently small
such that there is a fixed constant $\kappa\geq 1$ satisfying
\begin{equation}{\label{*clpeq74}}
\parallel {\cal P}\parallel\leq\kappa T/N.
\end{equation}

\item All the constraints of (CLP) and (DCLP) are satisfied for all $t\in {\cal P}$.

\item All the assumptions regarding the functions $a_{j}$, $c_{i}$, $B_{ij}$ and $K_{ij}$
are satisfied for all $t\in {\cal P}$.

\item All the discontinuities of $a_{j}$, $c_{i}$, $B_{ij}$ and $K_{ij}$
are contained in ${\cal P}$.

\item Remark~\ref{clpr42} is taken into account.
\end{itemize}

Given a partition ${\cal P}$ satisfying the above assumptions,
since each entry of ${\bf a}$ and ${\bf c}$ is piecewise continuous on $[0,T]$, i.e.,
each entry of ${\bf a}$ and ${\bf c}$ is continuous on each open interval
$(t_{u-1},t_{u})$ for $u=1,\cdots ,N$, we define
\begin{equation}{\label{dclp116}}
a_{j}^{(u)}=\inf_{t\in (t_{u-1},t_{u})}a_{j}(t)>-\infty\mbox{ and }
c_{i}^{(u)}=\inf_{t\in (t_{u-1},t_{u})}c_{i}(t)>-\infty.
\end{equation}
Then, by Remark~\ref{clpr42}, we see that
\begin{equation}{\label{extclp214}}
\tau\geq a_{j}^{(u)}\mbox{ and }\zeta\geq c_{i}^{(u)}.
\end{equation}
Now, we define the vectors ${\bf a}^{(u)}$ and ${\bf c}^{(u)}$ that are consisting
of $a_{j}^{(u)}$ and $c_{i}^{(u)}$ for $j=1,\cdots ,q$ and $i=1,\cdots ,p$,
respectively.

Since each entry of the time-dependent matrices $B$ and $K$ is piecewise continuous
on $[0,T]$ and $[0,T]\times [0,T]$, respectively, for $u=1,\cdots ,N$, we define
\begin{equation}{\label{dclp103}}
B_{ij}^{(u)}=\sup_{t\in (t_{u-1},t_{u})}B_{ij}(t)<+\infty\mbox{ and }
K_{ij}^{(u,v)}=\inf_{(t,s)\in (t_{u-1},t_{u})\times (t_{v-1},t_{v})}
K_{ij}(t,s)>-\infty .
\end{equation}
Then, by (\ref{*clpeq2}), (\ref{clpeq50}) and Remark~\ref{clpr42}, we see that
\begin{equation}{\label{extclpeq213}}
K_{ij}^{(u,v)}\leq\eta\mbox{ and }\sum_{i=1}^{p}K_{ij}^{(u,v)}\leq\nu
\end{equation}
We also define the matrices $B^{(u)}$ and $K^{(u,v)}$ that are consisting
of $B_{ij}^{(u)}$ and $K_{ij}^{(u,v)}$ for $j=1,\cdots ,q$ and $i=1,\cdots ,p$,
respectively.

Let ${\bf z}^{(u)}$ and ${\bf w}^{(u)}$ be the
$q$-dimensional and $p$-dimensional vectors, respectively.
We consider the following finite-dimensional linear programming problems:
\begin{eqnarray}
(\mbox{LP}^{(N)}) & \max & \sum_{u=1}^{N}\left (t_{u}-t_{u-1}\right )
({\bf a}^{(u)})^{\top}{\bf z}^{(u)}\nonumber\\
& \mbox{subject to} & B^{(1)}{\bf z}^{(1)}\leq
{\bf c}^{(1)}\nonumber\\
&& B^{(u)}{\bf z}^{(u)}\leq {\bf c}^{(u)}
+\sum_{v=1}^{u-1}\left (t_{v}-t_{v-1}\right )K^{(u,v)}
{\bf z}^{(v)}\mbox{ for }u=2,\cdots ,N\label{clpeq41}\\
&& {\bf z}^{(u)}\geq {\bf 0}\mbox{ for }u=1,\cdots ,N\nonumber.
\end{eqnarray}
and
\begin{eqnarray*}
(\mbox{DLP}^{(*N)}) & \min & \sum_{u=1}^{N}
({\bf c}^{(u)})^{\top}\widehat{\bf w}^{(u)}\\
& \mbox{subject to} & (B^{(u)})^{\top}\widehat{\bf w}^{(u)}
\geq \left (t_{u}-t_{u-1}\right ){\bf a}^{(u)}\\
&& \hspace{5mm}+\left (t_{u}-t_{u-1}\right )\sum_{v=u+1}^{N}(K^{(v,u)})^{\top}
\widehat{\bf w}^{(v)}\mbox{ for }u=1,\cdots ,N-1\\
&& (B^{(N)})^{\top}\widehat{\bf w}^{(N)}\geq
\left (t_{N}-t_{N-1}\right ){\bf a}^{(N)}\nonumber\\
&& \widehat{\bf w}^{(u)}\geq {\bf 0}
\mbox{ for }u=1,\cdots ,N.\nonumber
\end{eqnarray*}
Based on the following matrices
\[A{\bf z}={\small \left [\begin{array}{ccccc}
B^{(1)} & {\bf 0} & {\bf 0} & {\bf 0} & {\bf 0}\\
-(t_{1}-t_{0})K^{(2,1)} & B^{(2)} & & {\bf 0} & {\bf 0}\\
-(t_{1}-t_{0})K^{(3,1)} & -(t_{2}-t_{1})K^{(3,2)} & B^{(3)} & {\bf 0} & {\bf 0}\\
\vdots & \vdots & \vdots & \vdots & \vdots\\
-(t_{1}-t_{0})K^{(N,1)} & -(t_{2}-t_{1})K^{(N,2)}
& \cdots & -(t_{N}-t_{N-1})K^{(N,N-1)} & B^{(N)}
\end{array}\right ]
\left [\begin{array}{c}
{\bf z}^{(1)}\\
{\bf z}^{(2)}\\
{\bf z}^{(3)}\\
\vdots\\
{\bf z}^{(N)}
\end{array}\right ]}\]
and
\[A^{\top}\widehat{\bf w}=
{\small \left [\begin{array}{ccccc}
B^{(1)} & -(t_{1}-t_{0})K^{(2,1)} & -(t_{1}-t_{0})K^{(3,1)}
& \cdots & -(t_{1}-t_{0})K^{(N,1)}\\
{\bf 0} & B^{(2)} & -(t_{2}-t_{1})K^{(3,2)} & \cdots & -(t_{2}-t_{1})K^{(N,2)}\\
{\bf 0} & {\bf 0} & B^{(3)} & \cdots & -(t_{3}-t_{2})K^{(N,3)}\\
\vdots & \vdots & \vdots & \vdots & \vdots\\
{\bf 0} & {\bf 0} & \cdots & {\bf 0} & -(t_{N}-t_{N-1})K^{(N,N-1)}\\
{\bf 0} & {\bf 0} & \cdots & {\bf 0} & B^{(N)}
\end{array}\right ]
\left [\begin{array}{c}
\widehat{\bf w}^{(1)}\\
\widehat{\bf w}^{(2)}\\
\widehat{\bf w}^{(3)}\\
\vdots\\
\widehat{\bf w}^{(N-1)}\\
\widehat{\bf w}^{(N)}
\end{array}\right ],}\]
we see that $(\mbox{LP}^{(N)})$ and $(\mbox{DLP}^{(*N)})$ are finite-dimensional
primal and dual pair of linear programming problems. Now, let
\[{\bf w}^{(u)}=\left (\frac{1}{t_{u}-t_{u-1}}\right )\widehat{\bf w}^{(u)}.\]
Then, by dividing $t_{u}-t_{u-1}$ on both sides of the constraints of
dual problem $(\mbox{DLP}^{(*N)})$, we obtain the following equivalent problem
\begin{eqnarray}
(\mbox{DLP}^{(N)}) & \min & \sum_{u=1}^{N}\left (t_{u}-t_{u-1}\right )
({\bf c}^{(u)})^{\top}{\bf w}^{(u)}\nonumber\\
& \mbox{subject to} & (B^{(u)})^{\top}{\bf w}^{(u)}
\geq {\bf a}^{(u)}\nonumber\\
&& \quad\quad +\sum_{v=u+1}^{N}\left (t_{v}-t_{v-1}\right )(K^{(v,u)})^{\top}
{\bf w}^{(v)}
\mbox{ for }u=1,\cdots ,N-1\label{clpeq888}\\
&& (B^{(N)})^{\top}{\bf w}^{(N)}\geq {\bf a}^{(N)}\nonumber\\
&& {\bf w}^{(u)}\geq {\bf 0}\mbox{ for }u=1,\cdots ,N.\nonumber
\end{eqnarray}
In the sequel, we shall use this equivalent dual problem $(\mbox{DLP}^{(N)})$.

\begin{Pro}{\label{optp172}}
Suppose that each entry of ${\bf a}$, ${\bf c}$, $B$ and $K$ is piecewise continuous
on $[0,T]$ and $[0,T]\times [0,T]$, respectively. The following statements hold true.
\begin{enumerate}
\item [{\em (i)}] If ${\bf c}(t)\geq {\bf 0}$ a.e. in $[0,T]$,
then the primal problem $(\mbox{\em LP}^{(N)})$ is feasible.

\item [{\em (ii)}] Suppose that $\sum_{i=1}^{p}B_{ij}(t)>0$ a.e. in $[0,T]$ for each
$j=1,\cdots ,q$, and that there exists a constant $\sigma >0$ such that,
for each $i=1,\cdots ,p$ and
$j=1,\cdots ,q$, the following statement holds true a.e. in $[0,T]$:
\begin{equation}{\label{dclp118}}
B_{ij}(t)\neq 0\mbox{ implies }B_{ij}(t)\geq\sigma .
\end{equation}
Given a partition ${\cal P}=\{t_{0},t_{1},\cdots ,t_{N}\}$ of $[0,T]$, let
\begin{equation}{\label{clp2eq4}}
\mathfrak{w}_{u}=\frac{\tau}{\sigma}\cdot\left (1+\parallel {\cal P}\parallel\cdot
\frac{\nu}{\sigma}\right )^{N-u}\mbox{ for }u=1,\cdots ,N.
\end{equation}
We define the vector ${\bf w}^{(u)}$ with all entries $\mathfrak{w}_{u}$
for $u=1,\cdots ,N$. Then, $({\bf w}^{(1)},\cdots ,{\bf w}^{(N)})$ is a feasible solution of
problem $(\mbox{\em DLP}^{(N)})$; that is,
the dual problem $(\mbox{\em DLP}^{(N)})$ is feasible.
In other words, the strong duality theorem holds true between problems
$(\mbox{\em LP}^{(N)})$ and $(\mbox{\em DLP}^{(N)})$.
\end{enumerate}
\end{Pro}
\begin{Proof}
To prove part (i), since each entry of ${\bf c}$ is piecewise continuous on $[0,T]$,
according to the construction of partition ${\cal P}$, each entry $c_{i}$ is
continuous on the open subinterval $(t_{u-1},t_{u})$ for $u=1,\cdots ,N$.
Since $c_{i}(t)\geq 0$ a.e. on $(t_{u-1},t_{u})$, it follows that
$c_{i}(t)\geq 0$ for all $t\in (t_{u-1},t_{u})$ by the continuity.
From (\ref{dclp116}), we see that ${\bf c}^{(u)}\geq {\bf 0}$ for all $u=1,\cdots ,N$.
It is obvious that the primal problem $(\mbox{LP}^{(N)})$ is feasible with the
trivial feasible solution ${\bf z}^{(u)}={\bf 0}$ for $u=1,\cdots ,N$.

To prove part (ii), for each $j=1,\cdots ,q$, since the measure of open interval
$(t_{u-1},t_{u})$ is not zero, using the assumption for $B$ and referring to (\ref{dclp103}),
there exists $t^{*}\in (t_{u-1},t_{u})$ such that
$0\leq B_{ij}(t^{*})\leq B_{ij}^{(u)}$ for each $i=1,\cdots ,p$,
$\sum_{i=1}^{p}B_{ij}(t^{*})>0$ and the statement (\ref{dclp118})
is satisfied at $t^{*}$. Therefore,
there exists $i_{j}\in\{1,2,\cdots ,n\}$ such that $B_{i_{j}j}(t^{*})>0$, which implies
\[B_{i_{j}j}^{(u)}\geq B_{i_{j}j}(t^{*})\geq\sigma >0.\]
Since $B_{ij}^{(u)}\geq 0$ and $w_{i}^{(u)}=\mathfrak{w}_{u}\geq 0$ for
$i=1,\cdots ,p$ and $u=1,\cdots ,N$, we have
\begin{equation}{\label{clp2eq3}}
\sum_{i=1}^{p}B_{ij}^{(u)}\cdot w_{i}^{(u)}\geq B_{i_{j}j}^{(u)}\cdot
w_{i_{j}}^{(u)}=B_{i_{j}j}^{(u)}\cdot\frac{\tau}{\sigma}\cdot
\left (1+\parallel {\cal P}\parallel\cdot\frac{\nu}{\sigma}\right )^{N-u}
\geq\tau\cdot\left (1+\parallel {\cal P}\parallel\cdot\frac{\nu}{\sigma}\right )^{N-u}.
\end{equation}
Since
\begin{align*}
\sum_{i=1}^{p}\left (t_{v}-t_{v-1}\right )\cdot K_{ij}^{(v,u)}\cdot w_{i}^{(v)}
& \leq\sum_{i=1}^{p}\parallel {\cal P}\parallel\cdot K_{ij}^{(v,u)}\cdot
\frac{\tau}{\sigma}\cdot\left (1+\parallel {\cal P}\parallel\cdot
\frac{\nu}{\sigma}\right )^{N-v}\\
& \leq\parallel {\cal P}\parallel\cdot\frac{\nu\cdot\tau}{\sigma}
\left (1+\parallel {\cal P}\parallel\cdot\frac{\nu}{\sigma}\right )^{N-v}
\mbox{ (by (\ref{extclpeq213}))},
\end{align*}
it follows that, for $u=1,\cdots ,N-1$,
\begin{align}
& a_{j}^{(u)}+\sum_{i=1}^{p}\left [\sum_{v=u+1}^{N}
\left (t_{v}-t_{v-1}\right )\cdot K_{ij}^{(v,u)}\cdot w_{i}^{(v)}\right ]
=a_{j}^{(u)}+\sum_{v=u+1}^{N}\sum_{i=1}^{p}\left (t_{v}-t_{v-1}\right )\cdot
K_{ij}^{(v,u)}\cdot w_{i}^{(v)}\nonumber\\
& \quad\leq\tau+\sum_{v=u+1}^{N}\parallel {\cal P}\parallel\cdot
\frac{\nu\cdot\tau}{\sigma}\left (
1+\parallel {\cal P}\parallel\cdot\frac{\nu}{\sigma}\right )^{N-v}\nonumber\\
& \quad =\tau\cdot\left [1+\sum_{v=u+1}^{N}\parallel {\cal P}\parallel
\cdot\frac{\nu}{\sigma}\cdot\left (1+\parallel {\cal P}\parallel
\cdot\frac{\nu}{\sigma}\right )^{N-v}\right ]
=\tau\cdot\left (1+\parallel {\cal P}\parallel\cdot
\frac{\nu}{\sigma}\right )^{N-u}.\label{clp2eq5}
\end{align}
Therefore, from (\ref{clp2eq3}) and (\ref{extclp214}), we obtain
\[\sum_{i=1}^{p}B_{ij}^{(u)}\cdot w_{i}^{(u)}\geq a_{j}^{(u)}
+\sum_{i=1}^{p}\sum_{v=u+1}^{N}\left (t_{v}-t_{v-1}\right )\cdot
K_{ij}^{(v,u)}\cdot w_{i}^{(v)}\mbox{ for $u=1,\cdots ,N-1$}\]
and
\[\sum_{i=1}^{p}B_{ij}^{(N)}\cdot w_{i}^{(N)}\geq\tau\geq a_{j}^{(N)}.\]
This completes the proof.
\end{Proof}

\begin{Lem}{\label{*clpl72}}
Suppose that the set $\{x_{1},\cdots ,x_{N}\}$ satisfies
\[x_{1}\leq\theta_{1}\mbox{ and }
x_{u}\leq\theta_{1}+\theta_{2}\cdot\sum_{i=1}^{u-1}x_{i}\mbox{ for }u=2,\cdots ,N.\]
Then $x_{u}\leq\theta_{1}\cdot\left (1+\theta_{2}\right )^{u-1}$ for $u=1,\cdots ,N$.
\end{Lem}
\begin{Proof}
We are going to prove it by induction. For $u=1$, it is obviously true.
Suppose that $x_{u}\leq\theta_{1}\left (1+\theta_{2}\right )^{u-1}$
hold true for $u=2,\cdots ,N-1$. Then, we have
\[\sum_{i=1}^{N-1}x_{i}\leq\theta_{1}\cdot\sum_{i=1}^{N-1}
\left (1+\theta_{2}\right )^{i-1}=\frac{\theta_{1}
\left [(1+\theta_{2})^{N-1}-1\right ]}{\theta_{2}}.\]
Therefore, we obtain
\[x_{N}\leq\theta_{1}+\theta_{2}\sum_{i=1}^{N-1}x_{i}
\leq\theta_{1}+\theta_{2}\cdot\frac{\theta_{1}
\left [(1+\theta_{2})^{N-1}-1\right ]}{\theta_{2}}
=\theta_{1}\left (1+\theta_{2}\right )^{N-1}.\]
This completes the proof.
\end{Proof}

\begin{Pro}{\label{*clpp32}}
Suppose that each entry of ${\bf a}$, ${\bf c}$, $B$ and $K$ is piecewise continuous
on $[0,T]$ and $[0,T]\times [0,T]$, respectively,
and that there exist real-valued functions $\lambda_{i}$ satisfying
$0\leq\lambda_{i}(t)\leq 1$ a.e. in $[0,T]$ for $i=1,\cdots ,p$,
and a constant $\sigma >0$ such that
\begin{equation}{\label{dclp104}}
\min_{j=1,\cdots ,q}\left\{\sum_{i=1}^{p}\lambda_{i}(t)B_{ij}(t)
\right\}\geq\sigma\mbox{ a.e. in $[0,T]$}.
\end{equation}
If $({\bf z}^{(1)},\cdots ,{\bf z}^{(N)})$
is a feasible solution of primal problem $(\mbox{\em LP}^{(N)})$, then
\begin{equation}{\label{opteq164}}
\parallel{\bf z}^{(1)}\parallel\leq\frac{n\zeta}{\sigma}
\mbox{ and }\parallel{\bf z}^{(u)}\parallel
\leq\frac{n\zeta}{\sigma}\cdot\exp\left (\frac{n\nu\kappa T}{\sigma}\right )
\end{equation}
for $u=2,\cdots ,N$. In other words, the bound for the feasible solutions
is independent of the partition ${\cal P}$.
\end{Pro}
\begin{Proof}
For $u=1,\cdots ,N$, we define
\[\lambda_{i}^{(u)}=\mbox{ess}\sup_{t\in [t_{u-1},t_{u}]}\left |\lambda_{i}(t)\right |
=\inf\left\{k:\left |\lambda_{i}(t)\right |\leq k\mbox{ a.e. in }
[t_{u-1},t_{u}]\right\}.\]
Since $0\leq\lambda_{i}(t)\leq 1$ a.e. in $[0,T]$ for
$i=1,\cdots ,p$ by the assumption, it follows that
\begin{equation}{\label{dclp106}}
0\leq\lambda_{i}(t)\leq\lambda_{i}^{(u)}\mbox{ a.e. in }[t_{u-1},t_{u}]
\end{equation}
and $0\leq\lambda_{i}^{(u)}\leq 1$
for all $i=1,\cdots ,p$ and $u=1,\cdots ,N$. Since
$\sum_{j=1}^{q}B_{ij}^{(1)}z_{j}^{(1)}\leq c_{i}^{(1)}$ by the feasibility,
multiplying $\lambda_{i}^{(1)}\geq 0$ on both sides, we have
\begin{equation}{\label{dclp119}}
\sum_{i=1}^{p}\sum_{j=1}^{q}B_{ij}^{(1)}\lambda_{i}^{(1)}z_{j}^{(1)}\leq
\sum_{i=1}^{p}\lambda_{i}^{(1)}c_{i}^{(1)}.
\end{equation}
From (\ref{dclp103}), (\ref{dclp104}) and (\ref{dclp106}), we have
\begin{equation}{\label{dclp105}}
\min_{j=1,\cdots ,q}\left\{\sum_{i=1}^{p}\lambda_{i}^{(u)}B_{ij}^{(u)}
\right\}\geq\min_{j=1,\cdots ,q}\left\{\sum_{i=1}^{p}\lambda_{i}(t)B_{ij}(t)
\right\}\geq\sigma\mbox{ a.e. in $[t_{u-1},t_{u}]$}.
\end{equation}
Using (\ref{dclp119}), (\ref{dclp105}) and (\ref{extclp214}), we obtain
\[\sigma\cdot\parallel{\bf z}^{(1)}\parallel
=\sum_{j=1}^{q}\sigma\cdot z_{j}^{(1)}\leq\sum_{j=1}^{q}\left [z_{j}^{(1)}
\sum_{i=1}^{p}B_{ij}^{(1)}\lambda_{i}^{(1)}\right ]
\leq\sum_{i=1}^{p}\lambda_{i}^{(1)}c_{i}^{(1)}\leq\sum_{i=1}^{p}c_{i}^{(1)}
\leq n\zeta .\]
This shows that
\begin{equation}{\label{*clpeq73}}
\parallel{\bf z}^{(1)}\parallel\leq\frac{n\zeta}{\sigma}.
\end{equation}
From (\ref{clpeq41}), for each $u=2,\cdots ,N$ and $i=1,\cdots ,p$,
since $\lambda_{i}^{(u)}\geq 0$, we have
\begin{equation}{\label{copteq233}}
\sum_{i=1}^{p}\sum_{j=1}^{q}B_{ij}^{(u)}\lambda_{i}^{(u)}z_{j}^{(u)}\leq
\sum_{i=1}^{p}\lambda_{i}^{(u)}c_{i}^{(u)}+
\sum_{v=1}^{u-1}\sum_{i=1}^{p}\sum_{j=1}^{q}
\left (t_{v}-t_{v-1}\right )K_{ij}^{(u,v)}\lambda_{i}^{(u)}z_{j}^{(v)}.
\end{equation}
Since $0\leq\lambda_{i}^{(u)}\leq 1$, we obtain
\begin{align*}
\sigma\cdot\parallel{\bf z}^{(u)}\parallel
& =\sum_{j=1}^{q}\sigma\cdot z_{j}^{(u)}
\leq\sum_{j=1}^{q}\left [z_{j}^{(u)}
\sum_{i=1}^{p}B_{ij}^{(u)}\lambda_{i}^{(u)}\right ]\mbox{ (by (\ref{dclp105}))}\\
& \leq n\zeta +n\nu\parallel {\cal P}\parallel
\cdot\sum_{v=1}^{u-1}\parallel{\bf z}^{(v)}\parallel
\mbox{ (by (\ref{copteq233}), (\ref{clpeq50}) and (\ref{dclp103}))}.
\end{align*}
Let $\theta_{1}=n\zeta /\sigma$ and $\theta_{2}
=n\nu\parallel {\cal P}\parallel /\sigma$. We have
$\parallel{\bf z}^{(1)}\parallel\leq\theta_{1}$
by (\ref{*clpeq73}) and
\[\parallel{\bf z}^{(u)}\parallel\leq\theta_{1}
+\theta_{2}\cdot\sum_{v=1}^{u-1}\parallel{\bf z}^{(v)}\parallel
\mbox{ for $u=2,\cdots ,N$}.\]
According to Lemma~\ref{*clpl72}, we obtain
\begin{align*}
\parallel{\bf z}^{(u)}\parallel
& \leq\theta_{1}\left (1+\theta_{2}\right )^{u-1}
\leq\theta_{1}\left (1+\theta_{2}\right )^{N}\\
& \leq\theta_{1}\cdot\exp\left (\theta_{2}N\right )
\mbox{ (using the fact of $e^{t}\geq 1+t$)}\\
& =\frac{n\zeta}{\sigma}\cdot
\exp\left (\frac{n\nu N\parallel {\cal P}\parallel}{\sigma}\right )
\leq\frac{n\zeta}{\sigma}\cdot\exp\left (\frac{n\nu\kappa T}{\sigma}\right )
\mbox{ (using (\ref{*clpeq74}))}
\end{align*}
for $u=2,\cdots ,N$. This completes the proof.
\end{Proof}

\begin{Pro}{\label{coptp250}}
Suppose that each entry of ${\bf a}$, ${\bf c}$, $B$ and $K$ is piecewise continuous
on $[0,T]$ and $[0,T]\times [0,T]$, respectively,
and that there exists a constant $\sigma >0$ such that, for each $i=1,\cdots ,p$
and $j=1,\cdots ,q$, the following statement holds true a.e. in $[0,T]$:
\[B_{ij}(t)\neq 0\mbox{ implies }B_{ij}(t)\geq\sigma .\]
If $({\bf d}^{(1)},\cdots ,{\bf d}^{(N)})$ is a feasible
solution of dual problem $(\mbox{\em DLP}^{(N)})$, then there exists a
feasible solution $({\bf w}^{(1)},\cdots ,{\bf w}^{(N)})$ of dual problem
$(\mbox{\em DLP}^{(N)})$ such that
\begin{equation}{\label{opteq163}}
{\bf w}^{(u)}\leq {\bf d}^{(u)}\mbox{ and }
0\leq w_{i}^{(u)}\leq\frac{\tau}{\sigma}\cdot\exp
\left (\frac{\eta T}{\sigma}\right )
\end{equation}
for $u=1,\cdots ,N$ and $i=1,\cdots ,p$, where the bound of feasible solutions is
independent of the partition $\parallel {\cal P}\parallel$.
Moreover, if $({\bf d}^{(1)},\cdots ,{\bf d}^{(N)})$ is an optimal
solution of dual problem $(\mbox{\em DLP}^{(N)})$,
then $({\bf w}^{(1)},\cdots ,{\bf w}^{(N)})$ is also an
optimal solution of dual problem $(\mbox{\em DLP}^{(N)})$.
\end{Pro}
\begin{Proof}
Let
\begin{equation}{\label{*clpeq3}}
\rho (t)=\frac{\tau}{\sigma}\cdot\exp\left [\frac{\eta (T-t)}{\sigma}\right ].
\end{equation}
Then, we have
\begin{equation}{\label{*clpeq4}}
\sigma\rho (t)=\tau +\eta \cdot\int_{t}^{T}\rho (s)ds.
\end{equation}
From (\ref{*clpeq3}), for $u=1,\cdots ,N$, we define
\[\rho_{u}=\rho (t_{u})=\frac{\tau}{\sigma}\cdot\exp\left [
\frac{\eta (T-t_{u})}{\sigma}\right ].\]
By taking $t=t_{u}$ in (\ref{*clpeq4}), since $\rho$ is a decreasing function, we have
\begin{align}
\sigma\rho_{u} & =\tau +\eta\cdot\int_{t_{u}}^{T}\rho (s)ds
=\tau +\eta\cdot\sum_{v=u}^{N-1}\int_{t_{v}}^{t_{v+1}}\rho (s)ds
\geq\tau +\eta\cdot\sum_{v=u}^{N-1}\int_{t_{v}}^{t_{v+1}}\rho (t_{v+1})ds\nonumber\\
& \geq a_{j}^{(u)}+\sum_{v=u+1}^{N}(t_{v}-t_{v-1})K_{ij}^{(v,u)}\rho_{v}
\mbox{ (by (\ref{extclpeq213}) and (\ref{extclp214}))}.\label{copteq248}
\end{align}
According to constraint (\ref{clpeq888}), for $j=1,\cdots ,q$ and
$u=1,\cdots ,N-1$, we have
\begin{equation}{\label{copteq246}}
\sum_{i}B_{ij}^{(u)}d_{i}^{(u)}\geq a_{j}^{(u)}
+\sum_{v=u+1}^{N}(t_{v}-t_{v-1})K_{ij}^{(v,u)}d_{i}^{(v)}.
\end{equation}
For $u=1,\cdots ,N$, we define $w_{i}^{(u)}=\min\{d_{i}^{(u)},\rho_{u}\}$.
Since $K_{ij}^{(v,u)}\geq 0$, from (\ref{copteq246}), we obtain
\begin{equation}{\label{copteq247}}
\sum_{i}B_{ij}^{(u)}d_{i}^{(u)}\geq a_{j}^{(u)}
+\sum_{u=v+1}^{N}(t_{v}-t_{v-1})K_{ij}^{(v,u)}w_{i}^{(v)}.
\end{equation}
For each fixed $u$, we define the index sets $I_{\leq}=\{i:d_{i}^{(u)}\leq\rho_{u}\}$
and $I_{>}=\{i:d_{i}^{(u)}>\rho_{u}\}$ and consider
\[\sum_{i}B_{ij}^{(u)}w_{i}^{(u)}=\sum_{i\in I_{\leq}}
B_{ij}^{(u)}w_{i}^{(u)}+\sum_{i\in I_{>}}B_{ij}^{(u)}w_{i}^{(u)}.\]
Then, we have the following three cases.
\begin{itemize}
\item Suppose that $I_{>}=\emptyset$ (i.e., the second sum is zero).
Then, we see that $d_{i}^{(u)}=w_{i}^{(u)}$ for all $i$.
Therefore, from (\ref{copteq247}), we have
\[\sum_{i}B_{ij}^{(u)}w_{i}^{(u)}=\sum_{i}B_{ij}^{(u)}d_{i}^{(u)}\geq a_{j}^{(u)}
+\sum_{u=v+1}^{N}(t_{v}-t_{v-1})K_{ij}^{(v,u)}w_{i}^{(v)}.\]

\item Suppose that $I_{>}\neq\emptyset$ and $B_{ij}^{(u)}=0$ for all $i\in I_{>}$. Then
\begin{align*}
\sum_{i}B_{ij}^{(u)}w_{i}^{(u)}
& =\sum_{i\in I_{\leq}}B_{ij}^{(u)}d_{i}^{(u)}+\sum_{i\in I_{>}}B_{ij}^{(u)}\rho_{u}\\
& =\sum_{i\in I_{\leq}}B_{ij}^{(u)}d_{i}^{(u)}+\sum_{i\in I_{>}}B_{ij}^{(u)}d_{i}^{(u)}
=\sum_{i}B_{ij}^{(u)}d_{i}^{(u)}\\
& \geq a_{j}^{(u)}+\sum_{u=v+1}^{N}(t_{v}-t_{v-1})K_{ij}^{(v,u)}w_{i}^{(v)}
\mbox{ (by (\ref{copteq247}))}.
\end{align*}

\item Suppose that $I_{>}\neq\emptyset$, and that there exists $i^{*}\in I_{>}$ with
$B_{i^{*}j}^{(u)}\neq 0$. Then, by the definition of $B_{i^{*}j}^{(u)}$, there exists
$t^{*}\in (t_{u-1},t_{u})$ such that $B_{i^{*}j}(t^{*})\neq 0$. Since $B_{i^{*}j}(t)\geq 0$
a.e. in $[t_{u-1},t_{u}]$ and $B_{i^{*}j}$ is continuous on $(t_{u-1},t_{u})$,
we must have $B_{i^{*}j}(t)\geq 0$ for all $t\in (t_{u-1},t_{u})$.
The facts of $B_{i^{*}j}(t^{*})\neq 0$ and the continuity of $B_{i^{*}j}$ on $(t_{u-1},t_{u})$
imply that there exists a subset $T_{u}$ of
$(t_{u-1},t_{u})$ with nonzero measure such that $B_{i^{*}j}(t)>0$ on $T_{u}$.
By the assumption of $B$, we also have $B_{i^{*}j}(t)\geq\sigma$ a.e. in $T_{u}$.
Since $B_{i^{*}j}(t)\leq B_{i^{*}j}^{(u)}$ for all $t\in (t_{u-1},t_{u})$,
we conclude that $B_{i^{*}j}^{(u)}\geq\sigma$. Therefore, we obtain
\begin{align*}
\sum_{i}B_{ij}^{(u)}w_{i}^{(u)} & \geq\sum_{i\in I_{>}}B_{ij}^{(u)}w_{i}^{(u)}=
\sum_{i\in I_{>}}B_{ij}^{(u)}\rho_{u}\geq B_{i^{*}j}^{(u)}\rho_{u}\geq\sigma\rho_{u}\\
& \geq a_{j}^{(u)}+\sum_{u=v+1}^{N}(t_{v}-t_{v-1})K_{ij}^{(v,u)}w_{i}^{(v)}
\mbox{ (by (\ref{copteq248}) and $\rho_{v}\geq w_{i}^{(v)}\geq 0$)}.
\end{align*}
\end{itemize}
This shows that $({\bf w}^{(1)},\cdots ,{\bf w}^{(N)})$
is a feasible solution of dual problem $(\mbox{DLP}^{(N)})$.
Since $w_{i}^{(u)}\leq\rho_{u}$ for $u=1,\cdots ,N$, we also obtain (\ref{opteq163}).
Finally, since ${\bf w}^{(u)}\leq {\bf d}^{(u)}$
for $v=1,\cdots ,N$, if $({\bf d}^{(1)},\cdots ,{\bf d}^{(N)})$ is an optimal
solution of problem $(\mbox{DLP}^{(N)})$, then, considering the objective values,
we have
\[\sum_{u=1}^{N}(t_{u}-t_{u-1})({\bf c}^{(u)})^{\top}{\bf d}^{(u)}
\leq\sum_{u=1}^{N}(t_{u}-t_{u-1})({\bf c}^{(u)})^{\top}{\bf w}^{(u)}
\leq\sum_{u=1}^{N}(t_{u}-t_{u-1})({\bf c}^{(u)})^{\top}{\bf d}^{(u)},\]
which says that $({\bf w}^{(1)},\cdots ,{\bf w}^{(N)})$
is also an optimal solution of problem $(\mbox{DLP}^{(N)})$.
This completes the proof.
\end{Proof}

\section{Strong Duality Theorem}

Let ${\bf z}^{(u)}=(z_{1}^{(u)},\cdots ,z_{q}^{(u)})$ for $u=1,\cdots ,N$ be
feasible solutions of primal problem $(\mbox{LP}^{(N)})$ with the corresponding
partition ${\cal P}=\{0=t_{0},t_{1},\cdots ,t_{N}=T\}$ of $[0,T]$.
We define the step function
$\widehat{\bf z}(t)=(\widehat{z}_{1}(t),\cdots ,\widehat{z}_{q}(t))$ by
\begin{equation}{\label{copteq244}}
\widehat{z}_{j}(t)=\left\{\begin{array}{ll}
z_{j}^{(u)} & \mbox{if }t_{u-1}\leq t<t_{u}\mbox{ and }u=1,\cdots ,N\\
z_{j}^{(N)} & \mbox{if }t=T
\end{array}\right .\mbox{ for }j=1,\cdots ,q.
\end{equation}
Let ${\bf w}^{(u)}=(w_{1}^{(u)},\cdots ,w_{p}^{(u)})$ for $u=1,\cdots ,N$
be feasible solutions of dual problem $(\mbox{DLP}^{(N)})$.
We similarly define the step function
$\widehat{\bf w}(t)=(\widehat{w}_{1}(t),\cdots ,\widehat{w}_{p}(t))$ by
\begin{equation}{\label{opteq184}}
\widehat{w}_{i}(t)=\left\{\begin{array}{ll}
w_{i}^{(u)} & \mbox{if }t_{u-1}\leq t<t_{u}\mbox{ and }u=1,\cdots ,N\\
w_{i}^{(N)} & \mbox{if }t=T
\end{array}\right .\mbox{ for }i=1,\cdots ,p.
\end{equation}
Next, we present some useful lemmas for further discussion.

\begin{Lem}{\label{extdclp516}}
Suppose that each entry of ${\bf a}$, ${\bf c}$, $B$ and $K$ is piecewise continuous on $[0,T]$
and $[0,T]\times [0,T]$, respectively.
Given any $\epsilon >0$, we can take a sufficiently small
$\parallel {\cal P}\parallel$ with ${\cal P}=\{t_{0},t_{1},\cdots ,t_{N}\}$ that satisfies
$(\ref{*clpeq74})$ such that the following statements hold:
\begin{itemize}
\item $a_{j}(t)-a_{j}^{(u)}<\epsilon$ for $j=1,\cdots ,q$,
$t_{u-1}<t<t_{u}$ and $u=1,\cdots ,N$;

\item $c_{i}(t)-c_{i}^{(u)}<\epsilon$ for $i=1,\cdots ,p$,
$t_{u-1}<t<t_{u}$ and $u=1,\cdots ,N$;

\item $B_{ij}^{(u)}-B_{ij}(t)<\epsilon$ for $i=1,\cdots ,p$, $j=1,\cdots ,q$,
$t_{u-1}<t<t_{u}$ and $u=1,\cdots ,N$;

\item for fixed $t$ with $t_{u-1}<t<t_{u}$ and $u=1,\cdots ,N-1$, we have
$K_{ij}(s,t)-K_{ij}^{(v,u)}<\epsilon$ for  $i=1,\cdots ,p$, $j=1,\cdots ,q$,
$t_{v-1}<s<t_{v}$ and $v=u+1,\cdots ,N$.
\end{itemize}
\end{Lem}
\begin{Proof}
According to the construction of partition ${\cal P}$, we see that $a_{j}$ is
continuous on the open interval $E_{u}=(t_{u-1},t_{u})$.
We define the compact interval
\[E_{um}=\left [t_{u-1}+\frac{1}{m},t_{u}-\frac{1}{m}\right ].\]
Then
\begin{equation}{\label{clptmpceq16}}
E_{u}=\bigcup_{m=1}^{\infty}E_{um}\mbox{ and }
E_{um_{1}}\subseteq E_{um_{2}}\mbox{ for }m_{2}>m_{1}
\end{equation}
Since $E_{um}\subset E_{u}$, it follows that $a_{j}$ is continuous
on each compact interval $E_{um}$, which also means that $a_{j}$ is
uniformly continuous on each compact interval $E_{um}$.
Therefore, given any $\epsilon >0$, there exists $\delta >0$ such that
$|t_{1}-t_{2}|<\delta$ implies
\begin{equation}{\label{clptmpceq44}}
\left |a_{j}(t_{1})-a_{j}(t_{2})\right |<\epsilon
\mbox{ for any $t_{1},t_{2}\in E_{um}$.}
\end{equation}
Since the length of $E_{u}$ is less than or equal to
$\parallel {\cal P}\parallel\leq \kappa T/N$ by (\ref{*clpeq74}),
we can consider a sufficiently large $N_{0}\in\mathbb{N}$ such that
$\kappa T/N_{0}<\delta$. In this case, each length of $E_{u}$ for $l=1,\cdots ,p$
is less than $\delta$. In other words, if $N\geq N_{0}$,
then (\ref{clptmpceq44}) is satisfied for any $t_{1},t_{2}\in E_{um}$.
We consider the following cases.
\begin{itemize}
\item Suppose that the infimum $a_{j}^{(u)}$ is attained at
$t_{u}^{(*)}\in E_{u}$. From (\ref{clptmpceq16}), there exists $m^{*}$ such that
$t_{u}^{(*)}\in E_{um^{*}}$. Now, given any $t\in E_{u}$, we see that
$t\in E_{um_{0}}$ for some $m_{0}$. Let $m=\max\{m_{0},m^{*}\}$. From
(\ref{clptmpceq16}), it follows that $t,t_{u}^{(*)}\in E_{um}$. Then, we have
\[\left |a_{j}(t)-a_{j}^{(u)}\right |
=\left |a_{j}(t)-a_{j}\left (t_{u}^{(*)}\right )\right |<\epsilon\]
since the length of $E_{um}$ is less than $\delta$,
where $\epsilon$ is independent of $t$ because of the uniform continuity.

\item Suppose that the infimum $a_{j}^{(u)}$ is not attained at any point in
$E_{u}$. Since $a_{j}$ is continuous on the open interval $E_{u}$,
it follows that the infimum $a_{j}^{(u)}$ is either the
righthand limit or lefthand limit given by
\[a_{j}^{(u)}=\lim_{t\rightarrow t_{u-1}+}a_{j}(t)\mbox{ or }
a_{j}^{(u)}=\lim_{t\rightarrow t_{u}-}a_{j}(t).\]
Therefore, for sufficiently large $N_{0}$, i.e., the open interval $E_{u}$
is sufficiently small such that its length is less than $\delta$, we have
\[\left |a_{j}(t)-a_{j}^{(u)}\right |<\epsilon\]
for all $t\in E_{u}$.
\end{itemize}
From the above two cases, since $a_{j}(t)\geq a_{j}^{(u)}$ for all
$t\in E_{u}$, we conclude that $a_{j}(t)-a_{j}^{(u)}<\epsilon$.
The remaining cases can be similarly obtained. This completes the proof.
\end{Proof}

\begin{Lem}{\label{optl185}}
Suppose that the following conditions are satisfied:
\begin{itemize}
\item Suppose that each entry of ${\bf a}$, ${\bf c}$, $B$ and $K$ is piecewise continuous
on $[0,T]$ and $[0,T]\times [0,T]$, respectively;

\item there exists a constant $\sigma >0$ such that, for each $i=1,\cdots ,p$
and $j=1,\cdots ,q$, the following statement holds true a.e. in $[0,T]$:
\[B_{ij}(t)\neq 0\mbox{ implies }B_{ij}(t)\geq\sigma .\]
\end{itemize}
Let $({\bf d}^{(1)},\cdots ,{\bf d}^{(N)})$ be a feasible solution
of dual problem $(\mbox{\em DLP}^{(N)})$.
Given any $\epsilon >0$, there exists a sufficiently small
$\parallel {\cal P}\parallel$ with ${\cal P}=\{t_{0},t_{1},\cdots ,t_{N}\}$
which depends on $\epsilon$ such that there exists another
feasible solution $({\bf w}^{(1)},\cdots ,{\bf w}^{(N)})$
of dual problem $(\mbox{\em DLP}^{(N)})$ satisfying ${\bf w}^{(u)}\leq {\bf d}^{(u)}$
for $u=1,\cdots ,N$ and
\begin{equation}{\label{*clpeq30}}
\int_{0}^{T} {\bf c}^{\top}(t)\widehat{\bf w}(t)dt
\leq\sum_{u=1}^{N} (t_{u}-t_{u-1})({\bf c}^{(u)})^{\top}
{\bf w}^{(u)}+\epsilon ,
\end{equation}
where the step function $\widehat{\bf w}(t)$ is defined in $(\ref{opteq184})$.
If $({\bf d}^{(1)},\cdots ,{\bf d}^{(N)})$ is an optimal solution
of dual problem $(\mbox{\em DLP}^{(N)})$, then
$({\bf w}^{(1)},\cdots ,{\bf w}^{(N)})$ can be taken as the optimal
solution of dual problem $(\mbox{\em DLP}^{(N)})$.
\end{Lem}
\begin{Proof}
The existence of feasible solution $({\bf w}^{(1)},\cdots ,{\bf w}^{(N)})$
can be guaranteed by Proposition~\ref{coptp250}
From Lemma~\ref{extdclp516}, given any $\bar{\epsilon}>0$,
we can take a sufficiently small $\parallel {\cal P}\parallel$ with
${\cal P}=\{t_{0},t_{1},\cdots ,t_{N}\}$ that satisfies (\ref{*clpeq74})
such that $c_{i}(t)-c_{i}^{(u)}<\bar{\epsilon}$ for $i=1,\cdots ,p$,
$t_{u-1}<t<t_{u}$ and $u=1,\cdots ,N$, which implies
\begin{equation}{\label{opteq177}}
c_{i}(t) w_{i}^{(u)}-c_{i}^{(u)}w_{i}^{(u)}
\leq\bar{\epsilon}w_{i}^{(u)}\mbox{ for $t_{u-1}<t<t_{u}$ and $u=1,\cdots ,N$}
\end{equation}
by the fact of $w_{i}^{(u)}\geq 0$.
Since the integral does not be affected by the endpoints,
taking integrations from (\ref{opteq177}), we obtain
\[\sum_{i=1}^{p}\int_{0}^{T}c_{i}(t)\widehat{w}_{i}(t)dt-
\sum_{i=1}^{p}\sum_{u=1}^{N}(t_{u}-t_{u-1})c_{i}^{(u)} w_{i}^{(u)}
\leq\parallel {\cal P}\parallel\bar{\epsilon}\cdot\sum_{i=1}^{p}
\sum_{u=1}^{N} w_{i}^{(u)}.\]
By the boundedness of $ w_{i}^{(u)}$ as shown in (\ref{opteq163}), we also have
\begin{align*}
\int_{0}^{T}{\bf c}^{\top}(t)\widehat{\bf w}(t)dt-
\sum_{u=1}^{N} (t_{u}-t_{u-1})({\bf c}^{(u)})^{\top}{\bf w}^{(u)}
& \leq Np\bar{\epsilon}\parallel {\cal P}\parallel\frac{\tau}{\sigma}\cdot\exp
\left (\frac{\eta T}{\sigma}\right )\\
& \leq \kappa Tp\bar{\epsilon}\cdot\frac{\tau}{\sigma}\cdot\exp
\left (\frac{\eta T}{\sigma}\right )\mbox{ (by (\ref{*clpeq74}))}
\end{align*}
which implies, given any $\epsilon >0$, there exists a sufficiently small
$\parallel {\cal P}\parallel$ with ${\cal P}=\{t_{0},t_{1},\cdots ,t_{N}\}$ such that
\[\int_{0}^{T} {\bf c}^{\top}(t)\widehat{\bf w}(t)dt
-\sum_{u=1}^{N} (t_{u}-t_{u-1})({\bf c}^{(u)})^{\top}
{\bf w}^{(u)}\leq\epsilon .\]
This completes the proof.
\end{Proof}

\begin{Lem}{\label{optl183}}
Suppose that the following conditions are satisfied:
\begin{itemize}
\item each entry of ${\bf a}$, ${\bf c}$, $B$ and $K$ is piecewise continuous on $[0,T]$
and $[0,T]\times [0,T]$, respectively;

\item $K(s,t)\geq {\bf 0}$ a.e. in $[0,T]\times [0,T]$;

\item $\sum_{i=1}^{p}B_{ij}(t)>0$ a.e. in $[0,T]$ for each $j=1,\cdots ,q$;

\item there exists a constant $\sigma >0$ such that,
for each $i=1,\cdots ,p$
and $j=1,\cdots ,q$, the following statement holds true a.e. in $[0,T]$:
\[B_{ij}(t)\neq 0\mbox{ implies }B_{ij}(t)\geq\sigma .\]
\end{itemize}
Let $({\bf d}^{(1)},\cdots {\bf d}^{(N)})$ be a feasible solution
of dual problem $(\mbox{\em DLP}^{(N)})$.
Given any $\epsilon >0$, there exists a sufficiently small
$\parallel {\cal P}\parallel$ which depends on $\epsilon$ such that there exists another
feasible solution $({\bf w}^{(1)},\cdots {\bf w}^{(N)})$
of dual problem $(\mbox{\em DLP}^{(N)})$ satisfying
${\bf w}^{(u)}\leq {\bf d}^{(u)}$ for $u=1,\cdots ,N$ and
\[\sum_{i=1}^{p}B_{ij}(t)\widehat{w}_{i}(t)+\epsilon\geq a_{j}(t)
+\sum_{i=1}^{p}\int_{t}^{T}K_{ij}(s,t)\widehat{w}_{i}(s)ds\]
for all $t\in [0,T]\setminus {\cal P}$ and for $j=1,\cdots ,q$,
where the step function $\widehat{\bf w}(t)$ is defined in $(\ref{opteq184})$.
If $({\bf d}^{(1)},\cdots {\bf d}^{(N)})$ is an optimal
solution of $(\mbox{\em DLP}^{(N)})$, then
$({\bf w}^{(1)},\cdots {\bf w}^{(N)})$ can be taken as the optimal
solution of $(\mbox{\em DLP}^{(N)})$.
\end{Lem}
\begin{Proof}
The existence of feasible solution $({\bf w}^{(1)},\cdots ,{\bf w}^{(N)})$
can be guaranteed by Proposition~\ref{coptp250}
From Lemma~\ref{extdclp516}, given any $\bar{\epsilon}>0$,
we can take a sufficiently small $\parallel {\cal P}\parallel <\bar{\epsilon}$ with
${\cal P}=\{t_{0},t_{1},\cdots ,t_{N}\}$ that satisfies (\ref{*clpeq74})
such that, for $u=1,\cdots ,N$,
\begin{equation}{\label{opteq166}}
a_{j}(t)-a_{j}^{(u)}<\bar{\epsilon}\mbox{ and }B_{ij}^{(u)}-B_{ij}(t)<\bar{\epsilon}
\end{equation}
for $t_{u-1}<t<t_{u}$. Therefore, for $t_{u-1}<t<t_{u}$,using (\ref{opteq163}), we have
\begin{equation}{\label{opteq167}}
\sum_{i=1}^{p}B_{ij}^{(u)}w_{i}^{(u)}-\sum_{i=1}^{p}B_{ij}(t)\widehat{w}_{i}(t)
\leq\sum_{i=1}^{p}\bar{\epsilon}w_{i}^{(u)}
\leq p\bar{\epsilon}\cdot\frac{\tau}{\sigma}\cdot\exp
\left (\frac{\eta T}{\sigma}\right )\equiv\epsilon_{1}
\end{equation}
Also, from Lemma~\ref{extdclp516} again, for $t_{u-1}<t<t_{u}$, $u=1,\cdots ,N-1$
and $v=u+1,\cdots ,N$, we have $K_{ij}(s,t)-K_{ij}^{(v,u)}<\bar{\epsilon}$, which implies
\begin{align}
& \int_{t_{v-1}}^{t_{v}}K_{ij}(s,t)\cdot\widehat{w}_{i}(s)ds
-\left (t_{v}-t_{v-1}\right )K_{ij}^{(v,u)}w_{i}^{(v)}\nonumber\\
& \quad =\int_{t_{v-1}}^{t_{v}}\left (K_{ij}(s,t)-K_{ij}^{(v,u)}\right )w_{i}^{(v)}ds
\leq\bar{\epsilon}\left (t_{v}-t_{v-1}\right )w_{i}^{(v)}.\label{dclp130}
\end{align}
By referring to (\ref{*clpeq2}), for $t_{u-1}<t<t_{u}$, we obtain
\begin{equation}{\label{dclp108}}
\int_{t}^{t_{u}}K_{ij}(s,t)w_{i}^{(u)}ds
\leq\eta\parallel {\cal P}\parallel w_{i}^{(u)}
\leq\eta\bar{\epsilon}w_{i}^{(u)}.
\end{equation}
Now, for $u=1,\cdots ,N-1$, $t_{u-1}<t<t_{u}$ and $i=1,\cdots ,p$,
using (\ref{dclp130}) and (\ref{dclp108}), we have
\begin{align*}
& \int_{t}^{T}K_{ij}(s,t)\cdot\widehat{w}_{i}(s)ds
-\sum_{v=u+1}^{N}\left (t_{v}-t_{v-1}\right )K_{ij}^{(v,u)}w_{i}^{(v)}\\
& \quad\leq\bar{\epsilon}w_{i}^{(u)}\cdot\left [\eta +\sum_{v=u+1}^{N}\left (t_{v}-t_{v-1}\right )
\right ]\leq\bar{\epsilon}w_{i}^{(u)}\cdot\left (\eta +p\cdot\parallel {\cal P}\parallel \right )
\end{align*}
which implies, by using (\ref{opteq163}) and (\ref{*clpeq74}),
\begin{align}
& \sum_{i=1}^{p}\int_{t}^{T}K_{ij}(s,t)\cdot\widehat{w}_{i}(s)ds
-\sum_{i=1}^{p}\sum_{v=u+1}^{N}\left (t_{v}-t_{v-1}\right )
K_{ij}^{(v,u)}w_{i}^{(v)}\nonumber\\
& \quad\leq p\bar{\epsilon}\frac{\tau}{\sigma}\cdot\exp
\left (\frac{\eta T}{\sigma}\right )\cdot\left (\eta +\kappa T\right )
\equiv\epsilon_{2}.\label{opteq168}
\end{align}
For $u=1,\cdots ,N-1$ and $t_{u-1}<t<t_{u}$, using (\ref{opteq166}), (\ref{opteq167}),
(\ref{opteq168}) and the feasibility of $({\bf w}^{(1)},\cdots {\bf w}^{(N)})$,
we can obtain
\[-\sum_{i=1}^{p}B_{ij}(t)\widehat{w}_{i}(t)+a_{j}(t)
+\sum_{i=1}^{p}\int_{t}^{T}K_{ij}(s,t)\widehat{w}_{i}(s)ds
\leq\bar{\epsilon}+\epsilon_{1}+\epsilon_{2},\]
which shows that, for $i=1,\cdots ,p$, given $\bar{\epsilon}_{1}>0$, there exists a
sufficiently small $\parallel {\cal P}\parallel$ such that
\begin{equation}{\label{dclp110}}
\sum_{i=1}^{p}B_{ij}(t)\widehat{w}_{i}(t)+\bar{\epsilon}_{1}\geq a_{j}(t)
+\sum_{i=1}^{p}\int_{t}^{T}K_{ij}(s,t)\widehat{w}_{i}(s)ds.
\end{equation}
For $t_{N-1}<t<T$, using the similar argument, we can show that,
given $\bar{\epsilon}_{2}>0$, there exists a
sufficiently small $\parallel {\cal P}\parallel$ such that
\begin{equation}{\label{dclp111}}
\sum_{i=1}^{p}B_{ij}(t)\widehat{w}_{i}(t)+\bar{\epsilon}_{2}\geq a_{j}(t)
+\sum_{i=1}^{p}\int_{t}^{T}K_{ij}(s,t)\widehat{w}_{i}(s)ds.
\end{equation}
From (\ref{dclp110}) and (\ref{dclp111}), we complete the proof.
\end{Proof}

Let $M(\epsilon )$ and $\widehat{M}(\epsilon )$ be the optimal
objective values of $(\mbox{CLP}_{\epsilon})$ and
$(\mbox{DCLP}_{\epsilon})$, respectively. Under the assumptions of
Theorems~\ref{optt120} and \ref{p63},
we see that there exist optimal solutions $\bar{\bf z}^{(\epsilon)}$ and
$\bar{\bf w}^{(\epsilon)}$ of problems $(\mbox{CLP}_{\epsilon})$
and $(\mbox{DCLP}_{\epsilon})$, respectively, such that
\[M(\epsilon )=\int_{0}^{T} {\bf a}^{\top}(t)\bar{\bf z}^{(\epsilon)}(t)dt
\mbox{ and }\widehat{M}(\epsilon )
=\int_{0}^{T} {\bf c}^{\top}(t)\bar{\bf w}^{(\epsilon)}(t)dt.\]
Also, by taking $\epsilon =0$ in Theorems~\ref{optt120} and \ref{p63},
there exist optimal solutions ${\bf z}^{*}$ and ${\bf w}^{*}$
of problems (CLP) and (DCLP), respectively, such that
${\bf z}^{*}$ and ${\bf w}^{*}$ satisfy the following inequalities:
\begin{equation}{\label{*clpeq98}}
z_{j}^{*}(t)\leq\frac{p\cdot\zeta}{\sigma}\cdot
\exp\left (\frac{p\cdot\phi\cdot T}{\sigma}\right )\mbox{ a.e. in $[0,T]$
for each $j=1,\cdots ,q$}
\end{equation}
from (\ref{clpeq337}) by taking $\epsilon =0$, and
\begin{equation}{\label{*clpeq95}}
w_{i}^{*}(t)\leq\frac{\tau}{\sigma}\cdot\exp
\left (\frac{\nu\cdot T}{\sigma}\right )\mbox{ for all $t\in [0,T]$
and for each $i=1,\cdots ,p$}
\end{equation}
from (\ref{*clpeq68}).
Let $M$ and $\widehat{M}$ be the optimal objective values of (CLP) and (DCLP),
respectively. Then, we see that
\[M(0)=M=\int_{0}^{T} {\bf a}^{\top}(t){\bf z}^{*}(t)dt
\mbox{ and }\widehat{M}(0)=\widehat{M}=\int_{0}^{T} {\bf c}^{\top}(t){\bf w}^{*}(t)dt.\]
We are going to show that the functions $M(\epsilon )$ and
$\widehat{M}(\epsilon )$ are right-continuous at $0$.

\begin{Pro}{\label{*clpp70}}
Suppose that the following conditions are satisfied:
\begin{itemize}
\item each entry of ${\bf a}$, ${\bf c}$, $B$ and $K$ is piecewise continuous in $[0,T]$
and $[0,T]\times [0,T]$, respectively;

\item ${\bf c}(t)\geq {\bf 0}$ a.e. in $[0,T]$;

\item $K(t,s)\geq {\bf 0}$ a.e. in $[0,T]\times [0,T]$;

\item $\sum_{i=1}^{p}B_{ij}(t)>0$ a.e. in $[0,T]$ for each $j=1,\cdots ,q$;

\item there exists a constant $\sigma >0$ such that, for each $i=1,\cdots ,p$ and
$j=1,\cdots ,q$, the following statement holds true a.e. in $[0,T]$:
\[B_{ij}(t)\neq 0\mbox{ implies }B_{ij}(t)\geq\sigma .\]
\end{itemize}
Then, we have the following results.
\begin{enumerate}
\item [{\em (i)}] The function $M(\epsilon )$ is nondecreasing
and right-continuous at $0$, i.e., $M(0+)=M(0)$, and
\begin{equation}{\label{clpeq338}}
\lim_{\epsilon\rightarrow 0+}\int_{0}^{T}
{\bf a}^{\top}(t)\bar{\bf z}^{(\epsilon)}(t)dt=\int_{0}^{T}
{\bf a}^{\top}(t){\bf z}^{*}(t)dt,
\end{equation}
where ${\bf z}^{*}$ is an optimal solution of {\em (CLP)}
such that ${\bf z}^{*}(t)\geq {\bf 0}$ for all $t\in [0,T]$
and the inequalities in $(\ref{*clpeq98})$ are satisfied.
Moreover, the following results hold.
\begin{itemize}
\item If ${\bf c}(t)\geq {\bf 0}$ for all $t\in [0,T]$ and,
for each fixed $t_{0}\in [0,T]$, $K(t_{0},s)\geq {\bf 0}$ a.e. in $[0,T]$,
Then, there exists a common optimal solution
${\bf z}^{*}$ of {\em (CLP)} and $(\mbox{\em CLP}^{*})$
such that both problems have the same
optimal objective values and ${\bf z}^{*}$ satisfies the inequalities
$(\ref{*clpeq98})$;

\item If the conditions regarding the time-dependent matrix $B(t)$ are satisfied for all
$t\in [0,T]$, then the inequalities in $(\ref{*clpeq98})$ are satisfied for all $t\in [0,T]$.
\end{itemize}

\item [{\em (ii)}] The function $\widehat{M}(\epsilon )$ is nonincreasing and
right-continuous at $0$, i.e., $\widehat{M}(0+)=\widehat{M}(0)$, and
\begin{equation}{\label{clpeq339}}
\lim_{\epsilon\rightarrow 0+}\int_{0}^{T}
{\bf c}^{\top}(t)\bar{\bf w}^{(\epsilon)}(t)dt=\int_{0}^{T}
{\bf c}^{\top}(t){\bf w}^{*}(t)dt,
\end{equation}
where ${\bf w}^{*}$ is the optimal solution of {\em (DCLP)}
such that ${\bf w}^{*}(t)\geq {\bf 0}$ for all $t\in [0,T]$
and the inequalities $(\ref{*clpeq95})$ are satisfied.
If we further assume that the conditions regarding the time-dependent matrix
$B(t)$ are satisfied for all $t\in [0,T]$,
and that the function $\sum_{i=1}^{p}K_{ij}$ is bounded by $\nu$ and the function $a_{j}$
is bounded by $\tau$ for each $j=1,\cdots ,q$,
then there exists a common optimal solution
${\bf w}^{*}$ of {\em (DCLP)} and $(\mbox{\em DCLP}^{*})$
such that both problems have the same optimal objective values
and the inequalities in $(\ref{*clpeq95})$ are satisfied.
\end{enumerate}
\end{Pro}
\begin{Proof}
We want to show that $M(\epsilon )$ is nondecreasing.
For $\epsilon_{1}<\epsilon_{2}$, there exist optimal solutions
$\bar{\bf z}_{\epsilon_{1}}$ and $\bar{\bf z}_{\epsilon_{2}}$ satisfying
\begin{equation}{\label{opteq190}}
B(t)\bar{\bf z}_{\epsilon_{2}}(t)\leq
{\bf c}(t)+\mbox{\boldmath $\epsilon$}_{2}
+\int_{0}^{t} K(t,s)\bar{\bf z}_{\epsilon_{2}}(s)ds\mbox{ for all $t\in [0,T]$}
\end{equation}
and
\begin{equation}{\label{opteq191}}
B(t)\bar{\bf z}_{\epsilon_{1}}(t)\leq {\bf c}(t)+\mbox{\boldmath $\epsilon$}_{1}
+\int_{0}^{t} K(t,s)\bar{\bf z}_{\epsilon_{1}}(s)ds\leq
{\bf c}(t)+\mbox{\boldmath $\epsilon$}_{2}
+\int_{0}^{t} K(t,s)\bar{\bf z}_{\epsilon_{1}}(s)ds\mbox{ for all $t\in [0,T]$}.
\end{equation}
From (\ref{opteq191}), we see that $\bar{\bf z}_{\epsilon_{1}}$
is a feasible solution of $(\mbox{CLP}_{\epsilon_{2}})$.
Therefore, we obtain $M(\epsilon_{1})\leq M(\epsilon_{2})$,
since $(\mbox{CLP}_{\epsilon_{2}})$ is a maximization problem.
This shows that $M(\epsilon )$ is indeed nondecreasing, i.e., $M(0+)$ exists.
We also have
\begin{equation}{\label{opteq195}}
M=M(0)\leq M(0+).
\end{equation}
We consider the sequence $\{\epsilon_{k}\}_{k=1}^{\infty}$ such that
$\epsilon_{k}\rightarrow 0+$ as $k\rightarrow\infty$.
Using (\ref{clpeq337}) in Proposition~\ref{optt120*}
by taking $\lambda_{i}(t)=1$ for $i=1,\cdots ,p$ and Remark~\ref{*clpr97},
we see that the sequence $\{\bar{\bf z}^{(\epsilon_{k})}\}$ is
uniformly essentially bounded. Using part (i) of Proposition~\ref{*clpp47},
there exists a subsequence $\{\bar{\bf z}^{(\epsilon_{k_{r})}}\}$ which weakly
converges to some feasible solution $\bar{\bf z}^{(0)}\in L_{q}^{2}[0,T]$ of
$(\mbox{CLP}_{0})=\mbox{(CLP)}$. Moreover, there exists a feasible solution
${\bf z}^{*}$ of $(\mbox{CLP}_{0})=\mbox{(CLP)}$ such that
${\bf z}^{*}(t)\geq {\bf 0}$ for all $t\in [0,T]$
and ${\bf z}^{*}(t)=\bar{\bf z}^{(0)}(t)$ a.e. in $[0,T]$.
Therefore, using the weak convergence, we have
\begin{equation}{\label{opteq194}}
\int_{0}^{T}{\bf a}^{\top}(t){\bf z}^{*}(t)dt
=\int_{0}^{T}{\bf a}^{\top}(t)\bar{\bf z}^{(0)}(t)dt
=\lim_{r\rightarrow\infty}\int_{0}^{T}
{\bf a}^{\top}(t)\bar{\bf z}^{(\epsilon_{k_{r}})}(t)dt
=\lim_{r\rightarrow\infty}M(\epsilon_{k_{r}})=M(0+).
\end{equation}
Since ${\bf z}^{*}$ is a feasible solution of
$(\mbox{CLP}_{0})=\mbox{(CLP)}$, we also have
\begin{equation}{\label{opteq193}}
\int_{0}^{T}{\bf a}^{\top}(t){\bf z}^{*}(t)dt\leq M=M(0).
\end{equation}
Therefore, according to (\ref{opteq193}) and (\ref{opteq194}),
we obtain $M(0+)\leq M(0)$, which implies $M(0+)=M(0)=M$ by (\ref{opteq195}).
This says that ${\bf z}^{*}$ is an optimal solution of (CLP) and
proves equality (\ref{clpeq338}). If we further assume that
${\bf c}(t)\geq {\bf 0}$ for all $t\in [0,T]$ and,
for each fixed $t_{0}\in [0,T]$, $K(t_{0},s)\geq {\bf 0}$ a.e. in $[0,T]$,
using part (ii) of Proposition~\ref{*clpp47}, we see that
${\bf z}^{*}$ is also a feasible solution of $(\mbox{CLP}^{*})$.
Therefore, we conclude that ${\bf z}^{*}$ is also an optimal solution
of $(\mbox{CLP}^{*})$, since the feasible set of
$(\mbox{CLP}^{*})$ is contained in the feasible set of (CLP).
On the other hand, if the assumption regarding
the time-dependent matrix $B(t)$ is satisfied for all $t\in [0,T]$,
Then, according to (\ref{clpeq337}) in Proposition~\ref{optt120*},
the inequalities in (\ref{*clpeq98})
are satisfied for all $t\in [0,T]$. This proves part (i).

To prove part (ii), we can similarly show that $\widehat{M}(\epsilon )$ is
nonincreasing and $\widehat{M}=\widehat{M}(0)\geq\widehat{M}(0+)$.
From (\ref{*clpeq68}) in Theorem~\ref{p63}, we have
\[0\leq\bar{w}_{i}^{(\epsilon)}(t)\leq\frac{\tau}{\sigma}\cdot\exp
\left [\frac{\nu\cdot (T-t)}{\sigma}\right ]\leq\frac{\tau}{\sigma}\cdot\exp
\left (\frac{\nu\cdot T}{\sigma}\right )
\mbox{ for all $t\in [0,T]$},\]
which says that the sequence $\{\bar{\bf w}^{(\epsilon_{k})}\}_{k=1}^{\infty}$ is
uniformly essentially bounded. Using Proposition~\ref{*clpp48*},
there exists a subsequence $\{\widehat{\bf w}^{(\epsilon_{k_{r}})}\}_{r=1}^{\infty}$
which weakly converges to some feasible solution $\widehat{\bf w}^{(0)}\in L_{p}^{2}[0,T]$ of
$(\mbox{DCLP}_{0})=\mbox{(DCLP)}$ such that $\widehat{\bf w}^{(\epsilon_{k_{r}})}(t)
\leq\bar{\bf w}^{(\epsilon_{k_{r}})}(t)$ a.e. in $[0,T]$.
Moreover, there exists a feasible solution
${\bf w}^{*}$ of $(\mbox{DCLP}_{0})=\mbox{(DCLP)}$ such that
${\bf w}^{*}(t)\geq {\bf 0}$ for all $t\in [0,T]$
and ${\bf w}^{*}(t)=\widehat{\bf w}^{(0)}(t)$ a.e. in $[0,T]$.
Therefore, using the weak convergence, we have
\begin{align*}
\int_{0}^{T}{\bf c}^{\top}(t){\bf w}^{*}(t)dt
& =\int_{0}^{T}{\bf c}^{\top}(t)\widehat{\bf w}^{(0)}(t)dt=
\lim_{r\rightarrow\infty}\int_{0}^{T}
{\bf c}^{\top}(t)\widehat{\bf w}^{(\epsilon_{k_{r}})}(t)dt\\
& \leq\lim_{r\rightarrow\infty}\int_{0}^{T}
{\bf c}^{\top}(t)\bar{\bf w}^{(\epsilon_{k_{r}})}(t)dt
=\lim_{r\rightarrow\infty}\widehat{M}(\epsilon_{k_{r}})=\widehat{M}(0+).
\end{align*}
Since ${\bf w}^{*}$ is a feasible solution of
$(\mbox{DCLP}_{0})=\mbox{(DCLP)}$, we also have
\[\int_{0}^{T}{\bf c}^{\top}(t){\bf w}^{*}(t)dt\geq\widehat{M}=\widehat{M}(0).\]
Therefore, we obtain $\widehat{M}(0+)\geq\widehat{M}(0)$, which also implies
\[\widehat{M}(0+)=\int_{0}^{T}{\bf c}^{\top}(t){\bf w}^{*}(t)dt
=\widehat{M}(0)=\widehat{M},\]
This shows that ${\bf w}^{*}$ is an optimal solution of (DCLP),
and proves the equality (\ref{clpeq339}).

We further assume that the conditions regarding the time-dependent matrix
$B(t)$ are satisfied for all $t\in [0,T]$,
and that the function $\sum_{i=1}^{p}K_{ij}$ is bounded by $\nu$ and the function $a_{j}$
is bounded by $\tau$ for each $j=1,\cdots ,q$. Then,
part (iii) of Proposition~\ref{*clpp48*} says that we can take
${\bf w}^{*}$ as a feasible solution of $(\mbox{DCLP}^{*})$.
Since the feasible set of $(\mbox{DCLP}^{*})$ is contained in the
feasible set of $(\mbox{DCLP})$, it follows that ${\bf w}^{*}$
is an optimal solution of problem $(\mbox{DCLP}^{*})$. This completes the proof.
\end{Proof}

Now, we are in a position to prove the strong duality theorem.

\begin{Thm}{\label{optt203}}
{\em (Strong Duality Theorem)}
Suppose that the following conditions are satisfied:
\begin{itemize}
\item each entry of ${\bf a}$, ${\bf c}$, $B$ and $K$ is piecewise continuous in $[0,T]$
and $[0,T]\times [0,T]$, respectively;

\item ${\bf c}(t)\geq {\bf 0}$ a.e. in $[0,T]$;

\item $K(t,s)\geq {\bf 0}$ a.e. in $[0,T]\times [0,T]$;

\item the time-dependent matrix $B(t)$ satisfies the following conditions:
\begin{itemize}
\item $\sum_{i=1}^{p}B_{ij}(t)>0$ a.e. in $[0,T]$ for each $j=1,\cdots ,q$;

\item there exists a constant $\sigma >0$ such that, for each $i=1,\cdots ,p$ and
$j=1,\cdots ,q$, the following statement holds true a.e. in $[0,T]$:
\[B_{ij}(t)\neq 0\mbox{ implies }B_{ij}(t)\geq\sigma .\]
\end{itemize}
\end{itemize}
Then, there exist optimal solutions ${\bf z}^{*}$ and ${\bf w}^{*}$ of problems
{\em (CLP)} and {\em (DCLP)}, respectively, such that
\[\int_{0}^{T} {\bf a}^{\top}(t){\bf z}^{*}(t)dt=
\int_{0}^{T} {\bf c}^{\top}(t){\bf w}^{*}(t)dt,\]
where ${\bf z}^{*}(t)\geq {\bf 0}$ and ${\bf w}^{*}(t)\geq {\bf 0}$
for all $t\in [0,T]$, and the inequalities in
$(\ref{*clpeq98})$ and $(\ref{*clpeq95})$ are all satisfied.
Moreover, the following results hold.
\begin{itemize}
\item If ${\bf c}(t)\geq {\bf 0}$ for all $t\in [0,T]$ and,
for each fixed $t_{0}\in [0,T]$, $K(t_{0},s)\geq {\bf 0}$ a.e. in $[0,T]$,
then there exists a common optimal solution
${\bf z}^{*}$ of problems {\em (CLP)} and $(\mbox{\em CLP}^{*})$
such that both problems have the same optimal objective value.

\item If we further assume that the conditions regarding the time-dependent matrix
$B(t)$ are satisfied for all $t\in [0,T]$,
and that the function $\sum_{i=1}^{p}K_{ij}$ is bounded by $\nu$ and the function $a_{j}$
is bounded by $\tau$ for each $j=1,\cdots ,q$, then there exists a common optimal solution
${\bf w}^{*}$ of {\em (DCLP)} and $(\mbox{\em DCLP}^{*})$
such that both problems have the same optimal objective value
and the inequalities $(\ref{*clpeq98})$ are satisfied for all $t\in [0,T]$.
\end{itemize}
\end{Thm}
\begin{Proof}
Since the primal and dual pair of linear programming problems $(\mbox{LP}^{(N)})$ and
$(\mbox{DLP}^{(N)})$ are feasible by Proposition~\ref{optp172}, the strong duality theorem
says that there exist optimal solutions $(\bar{{\bf z}}^{(1)},\cdots ,\bar{{\bf z}}^{(N)})$
and $(\bar{{\bf d}}^{(1)},\cdots ,\bar{{\bf d}}^{(N)})$ of problems $(\mbox{LP}^{(N)})$ and
$(\mbox{DLP}^{(N)})$, respectively, such that
\begin{equation}{\label{opteq*187}}
\sum_{u=1}^{N}\left (t_{u}-t_{u-1}\right )({\bf a}^{(u)})^{\top}\bar{{\bf z}}^{(u)}
=\sum_{u=1}^{N}\left (t_{u}-t_{u-1}\right )({\bf c}^{(u)})^{\top}\bar{{\bf d}}^{(u)}.
\end{equation}
Since the integral does not be affected by the endpoints, from (\ref{copteq244}),
it is not hard to obtain
\begin{align}
\int_{0}^{T}{\bf a}^{\top}(t)\widehat{\bf z}(t)dt
& =\sum_{j=1}^{q}\int_{0}^{T}a_{j}(t)\widehat{z}_{j}(t)dt\nonumber\\
& \geq\sum_{j=1}^{q}\sum_{u=1}^{N}(t_{u}-t_{u-1})a_{j}^{(u)}\bar{z}_{j}^{(u)}
=\sum_{u=1}^{N} (t_{u}-t_{u-1})({\bf a}^{(u)})^{\top}\bar{\bf z}^{(u)}\label{*clpeq33}.
\end{align}

Considering $t\in [0,T]\setminus {\cal P}$,
since $K(t,s)$ is continuous on the open rectangles
$(t_{u-1},t_{u})\times (t_{v-1},t_{v})$ and $K(t,s)\geq {\bf 0}$ a.e. in
$(t_{u-1},t_{u})\times (t_{v-1},t_{v})$ for $u=1,\cdots ,N$ and $v=1,\cdots ,N$,
it follows that $K(t,s)\geq {\bf 0}$ for all $(t,s)\in
(t_{u-1},t_{u})\times (t_{v-1},t_{v})$ for $u=1,\cdots ,N$ and $v=1,\cdots ,N$,
which implies $K_{ij}^{(u,v)}\geq 0$ for $u=1,\cdots ,N$ and $v=1,\cdots ,N$.
Since $z_{j}^{(v)}\geq 0$ and $K_{ij}(t,s)\geq K_{ij}^{(u,v)}$ for $t_{u-1}<t<t_{u}$
and $t_{v-1}<s<t_{v}$, we have $K_{ij}(t,s)z_{j}^{(v)}\geq K_{ij}^{(u,v)}z_{j}^{(v)}$.
Since the integral does not be affected by the endpoints, for $t_{u-1}<t<t_{u}$,
we obtain
\begin{equation}{\label{dclp125}}
-\sum_{j=1}^{q}\int_{0}^{t}K_{ij}(t,s)\widehat{z}_{j}(s)ds\leq
-\sum_{j=1}^{q}\sum_{v=1}^{u-1}\left (t_{v}-t_{v-1}\right )K_{ij}^{(u,v)}z_{j}^{(v)}.
\end{equation}
Since $c_{i}(t)\geq c_{i}^{(u)}$ and $B_{ij}(t)\leq B_{ij}^{(u)}$,
using (\ref{dclp125}) and the feasibility of $({\bf z}^{(1)},\cdots ,{\bf z}^{(N)})$,
after some calculations, we can obtain the following inequalities
\[\sum_{j=1}^{q}B_{ij}(t)\widehat{z}_{j}(t)\leq c_{i}(t)
+\sum_{j=1}^{q}\int_{0}^{t}K_{ij}(t,s)\widehat{z}_{j}(s)ds\]
for all $t\in [0,T]\setminus {\cal P}$ and for $i=1,\cdots ,p$, which implies, for $\epsilon >0$,
\begin{equation}{\label{*clpeq32}}
B(t)\widehat{\bf z}(t)\leq {\bf c}(t)+\mbox{\boldmath $\epsilon$}
+\int_{0}^{t} K(t,s)\widehat{\bf z}(s)ds\mbox{ for all $t\in [0,T]\setminus {\cal P}$}
\end{equation}

Considering the dual problem $(\mbox{DLP}^{(N)})$,
applying Lemmas~\ref{optl185} and \ref{optl183}, there exists
a sufficiently small $\parallel {\cal P}\parallel$ which depends on $\epsilon$ such that
$(\bar{{\bf w}}^{(1)},\cdots ,\bar{{\bf w}}^{(N)})$ is an optimal solution of
problem $(\mbox{DLP}^{(N)})$ and the following inequalities are satisfied:
\begin{equation}{\label{*clpeq36}}
B^{\top}(t)\widehat{\bf w}(t)+\mbox{\boldmath $\epsilon$}\geq {\bf a}(t)
+\int_{t}^{T} K^{\top}(s,t)\widehat{\bf w}(s)ds
\mbox{ for all $t\in [0,T]\setminus {\cal P}$}
\end{equation}
and
\begin{equation}{\label{*clpeq202}}
\int_{0}^{T} {\bf c}^{\top}(t)\widehat{\bf w}(t)dt
\leq\sum_{u=1}^{N}\left (t_{u}-t_{u-1}\right )({\bf c}^{(u)})^{\top}
{\bf w}^{(u)}+\epsilon .
\end{equation}
From (\ref{opteq*187}), we also have
\begin{equation}{\label{opteq187}}
\sum_{u=1}^{N}\left (t_{u}-t_{u-1}\right )({\bf a}^{(u)})^{\top}\bar{{\bf z}}^{(u)}
=\sum_{u=1}^{N}\left (t_{u}-t_{u-1}\right )({\bf c}^{(u)})^{\top}\bar{{\bf w}}^{(u)}.
\end{equation}
The inequalities (\ref{*clpeq32}) and (\ref{*clpeq36}) say that
$\widehat{\bf z}$ and $\widehat{\bf w}$ are feasible solutions
of problems $(\mbox{CLP}_{\epsilon})$ and $(\mbox{DCLP}_{\epsilon})$, respectively.
By Theorems~\ref{optt120} and \ref{p63},
we see that there exist optimal solutions $\bar{\bf z}^{(\epsilon)}$ and
$\bar{\bf w}^{(\epsilon)}$ of $(\mbox{CLP}_{\epsilon})$
and $(\mbox{DCLP}_{\epsilon})$, respectively, such that
\begin{equation}{\label{opteq181}}
\int_{0}^{T} {\bf a}^{\top}(t)\bar{\bf z}^{(\epsilon)}dt
\geq\int_{0}^{T} {\bf a}^{\top}(t)\widehat{\bf z}(t)dt
\mbox{ and }\int_{0}^{T} {\bf c}^{\top}(t)\bar{\bf w}^{(\epsilon)}dt
\leq\int_{0}^{T} {\bf c}^{\top}(t)\widehat{\bf w}(t)dt.
\end{equation}
Using (\ref{opteq181}), (\ref{*clpeq33}) and (\ref{*clpeq202}), we have
\begin{equation}{\label{opteq182}}
\int_{0}^{T} {\bf a}^{\top}(t)\bar{\bf z}^{(\epsilon)}dt
\geq\sum_{u=1}^{N}\left (t_{u}-t_{u-1}\right )({\bf a}^{(u)})^{\top}
\bar{{\bf z}}^{(u)}
\end{equation}
and
\begin{equation}{\label{opteq186}}
\int_{0}^{T} {\bf c}^{\top}(t)\bar{\bf w}^{(\epsilon)}(t)dt
\leq\sum_{u=1}^{N}\left (t_{u}-t_{u-1}\right )({\bf c}^{(u)})^{\top}
\bar{{\bf w}}^{(u)}+\epsilon .
\end{equation}
From (\ref{opteq187}), (\ref{opteq182}) and (\ref{opteq186}), we obtain
\begin{equation}{\label{opteq188}}
\int_{0}^{T} {\bf c}^{\top}(t)\bar{\bf w}^{(\epsilon)}(t)dt
\leq\int_{0}^{T} {\bf a}^{\top}(t)\bar{\bf z}^{(\epsilon)}dt+\epsilon .
\end{equation}
Since $M(\epsilon )\uparrow M(0)=M$ and
$\widehat{M}(\epsilon )\downarrow\widehat{M}(0)=\widehat{M}$
as $\epsilon\rightarrow 0+$ by Proposition~\ref{*clpp70},
for each $\epsilon >0$, we obtain
\begin{equation}{\label{opteq196}}
\int_{0}^{T}{\bf a}^{\top}(t)\bar{\bf z}^{(\epsilon)}(t)dt\leq
\int_{0}^{T}{\bf a}^{\top}(t){\bf z}^{*}(t)dt\mbox{ and }
\int_{0}^{T}{\bf c}^{\top}(t)\bar{\bf w}^{(\epsilon)}(t)dt\geq
\int_{0}^{T}{\bf c}^{\top}(t){\bf w}^{*}(t)dt.
\end{equation}
By (\ref{opteq188}) and (\ref{opteq196}), we have
\[\int_{0}^{T}{\bf c}^{\top}(t){\bf w}^{*}(t)dt\leq
\int_{0}^{T}{\bf a}^{\top}(t){\bf z}^{*}(t)dt+\epsilon ,\]
which implies
\[\int_{0}^{T}{\bf c}^{\top}(t){\bf w}^{*}(t)dt\leq
\int_{0}^{T}{\bf a}^{\top}(t){\bf z}^{*}(t)dt,\]
since $\epsilon$ can be any positive number.
By the weak duality Theorem~\ref{optt198}, we conclude that
\[\int_{0}^{T}{\bf c}^{\top}(t){\bf w}^{*}(t)dt=
\int_{0}^{T}{\bf a}^{\top}(t){\bf z}^{*}(t)dt.\]
Also, the inequalities regarding the bounds of optimal solutions
${\bf z}^{*}$ and ${\bf w}^{*}$ follow from Proposition~\ref{*clpp70} immediately.

If ${\bf c}(t)\geq {\bf 0}$ for all $t\in [0,T]$ and,
for each fixed $t_{0}\in [0,T]$, $K(t_{0},s)\geq {\bf 0}$ a.e. in $[0,T]$,
then part (i) of Proposition~\ref{*clpp70} says that ${\bf z}^{*}$ is a common optimal
solution of problems {(CLP)} and $(\mbox{CLP}^{*})$ such that both problems have the
same optimal objective value.

Finally, we further assume that the conditions regarding the time-dependent matrix
$B(t)$ are satisfied for all $t\in [0,T]$,
and that the function $\sum_{i=1}^{p}K_{ij}$ is bounded by $\nu$ and the function $a_{j}$
is bounded by $\tau$ for each $j=1,\cdots ,q$.
Then, part (ii) of Proposition~\ref{*clpp70} says that ${\bf w}^{*}$ is a common optimal
solution of problems {(DCLP)} and $(\mbox{DCLP}^{*})$ such that both problems have the
same optimal objective value. Also, part (i) of Proposition~\ref{*clpp70} says that
the inequalities $(\ref{*clpeq98})$ are satisfied for all $t\in [0,T]$.
This completes the proof.
\end{Proof}

\end{document}